\documentclass[12pt,a4paper]{article}
\usepackage{graphicx,epsfig}
\usepackage{amsfonts}
\usepackage{amssymb}
\usepackage{amsmath}

\usepackage{times}
\usepackage{amsmath,amsthm}
\usepackage{graphics}

\newtheorem{theo}{Theorem}

\newtheorem{cor}{Corollary}
\newtheorem{rem}{Remark}
\newtheorem{ex}{Example}
\newtheorem{defn}{Definition}
\def\RR{\mathbb R}

\def\pmatrix{ \left( \begin{array} }
\def\endpmatrix{ \end{array} \right) }
\def\sigmd{{\dot\sigma}}

\def\no{\noindent}

\def\d{{\rm d}}

\def\pmatrix{ \left( \begin{array} }
\def\endpmatrix{ \end{array} \right) }
\def\aa{\alpha}

\def\bfy{\boldsymbol{y}}

\def\d2dxx{\frac{\partial^2}{\partial x^2}}

\def\no{\noindent}
\def\diag{{\rm diag}}

\def\QED{~\mbox{$\Box$}}

\def\phi{\varphi}

\def\I{{\cal I}}
\def\P{{\cal P}}
\def\D{\Lambda}
\def\O{\Omega}

\title{Hamiltonian Boundary Value Methods \\
{\large (Energy Conserving Discrete Line Integral
Methods)\thanks{Work developed within the project ``Numerical
methods and software for differential \newline equations''.}}}

\author{Luigi Brugnano\thanks{Dipartimento di Matematica ``U.Dini",
Viale Morgagni 67/A , I-50134 Firenze, Italy. \newline
(luigi.brugnano@unifi.it).} \and Felice
Iavernaro\thanks{Dipartimento di Matematica, Via Orabona 4,
I-70125 Bari, Italy. (felix@dm.uniba.it).} \and Donato
Trigiante\thanks{Dipartimento di Energetica ``S.Stecco", Via
Lombroso 6/17, I-50134 Firenze, Italy.  \newline (trigiant@unifi.it).}}

\date{~\\ \small To apper in: {\em Journal of Numerical Analysis,
Industrial and Applied Mathematics\\ (JNAIAM)}.\\ Received: October 25, 2009.
~ Accepted (in revised form): April 15, 2010.}

\begin{document}

\maketitle

\begin{abstract}
Recently, a new family of integrators (Hamiltonian Boundary Value
Methods) has been introduced, which is able to precisely conserve
the energy function of polynomial Hamiltonian systems and to
provide a \textit{practical} conservation of the energy in the
non-polynomial case.

We settle the definition and the theory of such methods in a more
general framework.  Our aim is on the one hand to give account of
their good behavior when applied to general Hamiltonian systems
and, on the other hand, to find out what are  the \textit{optimal}
formulae, in relation to the choice of the polynomial basis and of
the distribution of the nodes. Such analysis is based upon the
notion of \textit{extended collocation conditions} and the
definition of \textit{discrete line integral}, and is carried out
by looking at the limit of such family of methods as the number of
the so called {\em silent stages} tends to infinity.

\medskip
\no{\bf Keywords:} Hamiltonian problems, exact conservation of the
Hamiltonian, energy conservation, Hamiltonian Boundary Value
Methods, HBVMs, discrete line integral.

\medskip
\no{\bf MSC} 65P10, 65L05.
\end{abstract}

\section{Introduction}

We consider canonical Hamiltonian problems in the form
\begin{equation}\label{hamilode}
\dot y = J\nabla H(y),  \qquad y(t_0) = y_0\in\RR^{2m},
\end{equation}
where $J$ is a skew-symmetric constant matrix, and the Hamiltonian
$H(y)$ is assumed to be sufficiently differentiable. For its
numerical integration, the problem is to find numerical methods
which preserve $H(y)$ along the discrete solution $\{y_n\}$, since
this property holds for the continuous solution $y(t)$.

So far, many attempts have been made inside the class of
Runge-Kutta methods, the most successful of them being that of
imposing the symplecticity of the discrete map, considering that,
for the continuous flow, symplecticity implies the conservation of
$H(y)$. Concerning symplectic integrators, a backward error
analysis permits to prove that they exactly conserve a modified
Hamiltonian, even though this fact clearly does not always
guarantee a proper qualitative behavior of the discrete orbits.

On the other hand, it is possible to follow different approaches
to derive geometric integrators which are energy-preserving. This
has been done, for example, in the pioneering work \cite{Gon96},
and later in \cite{McL99}, where {\em discrete gradient methods}
are introduced and studied. An additional example of
energy-preserving method is the {\em Averaged Vector Field (AVF)}
method defined in \cite{Qui08} (see also \cite{M2AN}). By the way,
the latter method can be retrieved by the methods here studied.

More recently, in  \cite{BIT09} a new family of one-step methods
has been introduced, capable of providing a numerical solution
$\{y_n\}$ of \eqref{hamilode}, along which the energy function
$H(y)$ is precisely conserved, in the case where this function is
a polynomial (see also \cite{IP1,IT2,HBVMHome}).

These methods, named {\em Hamiltonian Boundary Value Methods}
({\em HBVMs} hereafter), may be also thought of as Runge-Kutta
methods where the internal stages are split into two categories:
\begin{itemize}
\item[-] the \textit{fundamental stages}, whose number, say $s$,
is related to the order of the method;
\item[-] the \textit{silent stages}, whose number, say $r$, has
to be suitably selected in order to assure the energy conservation
property for a polynomial $H(y)$ of given degree $\nu$; the higher
is $\nu$, the higher must be $r$.
\end{itemize}
The resulting method is denoted by HBVM$(k,s)$,\footnote{The
denomination {\em HBVM} with $k$ {\em steps} and {\em degree} $s$
was used in \cite{BIT09}.} where $k=s+r$ is the total number of
{\em unknown} stages.

In \cite{BIT09,IT2} it has also been shown that these new methods
provide a practical conservation of the energy even in the
non-polynomial case: the term ``practical'' means that, in many
general situations, when the number of silent stages is high
enough, the method makes no distinction between the function
$H(y)$ and its polynomial approximation, being the latter in a
neighborhood of size $\varepsilon$ of the former, where
$\varepsilon$ denotes the machine precision.

Another relevant issue to be mentioned is that the computational
cost for the solution of the associated nonlinear system is
essentially independent of the number of silent stages, and only
depends on $s$ (see \cite{BIT09,HBVMHome}). This comes from the
fact that the silent stages are actually linear combinations of
the fundamental stages.

These two aspects motivate the following question: \textit{what
is,  if any, the limit method when the number of silent stages
grows to infinity?}

This question was first posed by Ernst Hairer,\footnote{During the
international conference ``ICNAAM 2009'',  Rethymno, Crete,
Greece, 18-22 September 2009, after the talks, by the first two
authors, where HBVMs were presented.} who also provided a partial
answer by stating formulae \eqref{hbvm_hairer}, which he called
\textit{Energy Preserving  variant of Collocation Methods} ({\em
EPCMs}, hereafter) \cite{hairer}. We provide a proof of his
statement by clarifying the connection between the limit formulae
and HBVMs: we show that actually one can define several different
limit methods,\footnote{In the sense that they generate different
discrete problems.} each one associated to the specific polynomial
basis, as well as to the choice of the abscissae distribution,
used to construct the sequence of HBVMs. For example,
EPCMs are based upon the use of Lagrange polynomials, while,
working with the shifted Legendre basis, yields to different limit
methods, that we have called {\em Infinity Hamiltonian Boundary
Value Methods} (in short, {\em $\infty$-HBVMs} or {\em
HBVM$(\infty,s)$}, being $s$ the number of the {\em unknown}
fundamental stages).

Our aim in this paper is threefold:
\begin{enumerate}
\item  We settle the definition  of HBVMs in a more general framework,
also deriving the general formulation of the limit formulae
$$
\lim_{k \rightarrow \infty}
\mbox{HBVM}(k,s).
$$
In particular, we show  that such limit coincides with EPCMs if
the Lagrange polynomial basis is used (Section \ref{sect:def}).

\item In Section \ref{sect:inf}, we introduce  the new class of $\infty$-HBVMs,
 which are the limit formulae corresponding to the HBVMs based upon the shifted
 Legendre polynomial basis. We prove that the order of such formulae is the same
 as the Gauss-Legendre methods, that is $2s$ (where $s$ is the number of the
 unknown fundamental stages).

\item  We  mention the case where $H(y)$ belongs to vector spaces different from
that of polynomials, thus providing a natural (and trivial)
generalization of the original formulae (see Section
\ref{sect:gen}).  Moreover, in the polynomial case, we determine
the {\em optimal} distribution of the nodes (Section
\ref{sect:gauss}).
\end{enumerate}

We stress that any finite approximation of EPCMs or $\infty$-HBVMs
based on quadratures leads back to  HBVM($k$,$s$) methods, for $k$
high enough.

We address all the points listed above, by slightly modifying the
approach followed to define the class of HBVMs in \cite{BIT09}.

\section{Reformulation of Hamiltonian BVMs}
\label{sect:def}

The key formula which HBVMs rely on, is the {\em line integral}
and the related property of conservative vector fields:
\begin{equation}\label{Hy}
H(y_1) - H(y_0) = h\int_0^1 \sigmd(t_0+\tau h)^T\nabla
H(\sigma(t_0+\tau h))\mathrm{d}\tau, \qquad \mbox{for any $y_1 \in
\RR^{2m}$},
\end{equation}
where $\sigma$ is any smooth function such that
\begin{equation}
\label{sigma}\sigma(t_0) = y_0, \qquad\sigma(t_0+h) =
y_1.
\end{equation}
Here we consider the case where $\sigma(t)$ is a polynomial (of
degree at most $s$), yielding an approximation to the true
solution $y(t)$ in the time interval $[t_0,t_0+h]$. The numerical
approximation for the subsequent time-step, $y_1$, is then defined
by (\ref{sigma}). After introducing a set of $s$ distinct
abscissae  $c_{1},\ldots ,c_{s}$,  ($0< c_{i}\le 1$),\footnote{As
a convention, when $c=0$ is to be considered, as in the case of
the Lobatto abscissae in $[0,1]$, then $c_0=0$ is formally added
to the abscissae $c_1,\dots,c_s$, and the subsequent formulae are modified
accordingly.} we set \begin{equation}\label{Yi}Y_i=\sigma(t_0+c_i h), \qquad
i=1,\dots,s,\end{equation}
so that $\sigma(t)$ may be thought
of as an interpolation polynomial, $Y_i$, $i=1,\dots,s$, being the
internal stages.

Let us consider the following expansions of
$\dot \sigma(t)$ and $\sigma(t)$ for $t\in [t_0,t_0+h]$:
\begin{equation}
\label{expan} \dot \sigma(t_0+\tau h) = \sum_{j=1}^{s} \gamma_j
P_j(\tau), \qquad \sigma(t_0+\tau h) = y_0 + h\sum_{j=1}^{s}
\gamma_j \int_{0}^\tau P_j(x)\,\mathrm{d}x,
\end{equation}
where $\{P_j(t)\}$ is any suitable basis of the vector space of
polynomials of degree at most $s-1$ and  the  (vector)
coefficients $\{\gamma_j\}$ are to be determined.\footnote{More
general function spaces will be considered in the sequel.}
Before proceeding, one important remark is in order.

\begin{rem}
\label{rem:basis} As will be clear in a while, we observe that the
numerical method which the following procedure will define is
``basis-dependent'', in that to different choices of the basis
$\{P_j(t)\}$ there will, in general, correspond different
numerical methods. In this section, in order to let the theory be
presented as general as possible, we leave the basis not better
specified. This will allow us to achieve the results listed at
point 1. in the introduction. The question about how to choose the
basis properly is faced in Section \ref{sect:inf}, where
$\infty$-HBVMs will be introduced. Therefore, just in the present
section, to avoid confusion, we will always specify what is the
basis we are working with. This will be not necessary anymore
starting from Section \ref{sect:inf}, after determining the
\textit{optimal} basis.
\end{rem}

In this section we assume that $H(y)$ is a polynomial, which
implies that the integrand in \eqref{Hy} is also a polynomial so
that the line integral can be exactly computed by means of a
suitable quadrature formula. It is easy to observe that in
general, due to the high degree of the integrand function,  such
quadrature formula cannot be solely based upon the available
abscissae $\{c_i\}$: one needs to introduce an additional set of
abscissae, $\hat c_1, \dots,\hat c_r$, distinct from the nodes
$\{c_i\}$, in order to make the quadrature formula
exact:

\begin{eqnarray} \label{discr_lin}
\displaystyle \lefteqn{\int_0^1 \sigmd(t_0+\tau h)^T\nabla
H(\sigma(t_0+\tau h))\mathrm{d}\tau   =}\\ &&
\sum_{i=1}^s \beta_i \sigmd(t_0+c_i h)^T\nabla H(\sigma(t_0+c_i
h)) + \sum_{i=1}^r \hat \beta_i \sigmd(t_0+\hat c_i h)^T\nabla
H(\sigma(t_0+\hat c_i h)), \nonumber
\end{eqnarray}
where $\beta_i$, $i=1,\dots,s$, and $\hat \beta_i$, $i=1,\dots,r$,
denote the weights of the quadrature formula corresponding to the
abscissae $\{c_i\}$ and $\{\hat c_i\}$, respectively, i.e.,
\begin{eqnarray}\nonumber
\beta_i &=& \int_0^1\left(\prod_{ j=1,j\ne i}^s
\frac{t-c_j}{c_i-c_j}\right)\left(\prod_{j=1}^r
\frac{t-\hat c_j}{c_i-\hat c_j}\right)\mathrm{d}t, \qquad i = 1,\dots,s,\\
\label{betai}\\ \nonumber \hat\beta_i &=& \int_0^1\left(\prod_{
j=1}^s \frac{t-c_j}{\hat c_i-c_j}\right)\left(\prod_{ j=1,j\ne
i}^r \frac{t-\hat c_j}{\hat c_i-\hat c_j}\right)\mathrm{d}t,
\qquad i = 1,\dots,r.
\end{eqnarray}

According to \cite{IT2}, the right-hand side of \eqref{discr_lin}
is called \textit{discrete line integral}, while the vectors
\begin{equation}\label{hYi}
\hat Y_i = \sigma(t_0+\hat c_i h), \qquad i=1,\dots,r,
\end{equation} are called \textit{silent stages}: they just serve to
increase, as much as one likes, the degree of precision of the
quadrature formula, but they are not to be regarded as unknowns
since, from \eqref{expan}, they can be expressed in terms of
linear combinations of the \textit{fundamental stages} (\ref{Yi}).

In \cite{BIT09}, the method HBVM($k$,$s$), with $k=s+r$ is then
defined by  substituting the quantities in \eqref{expan} into the
right-hand side of \eqref{discr_lin} and by choosing the unknowns
$\{\gamma_j\}$ in order that the resulting expression vanishes.

Instead of carrying out our computation on the right-hand side of
\eqref{discr_lin}, as was done in \cite{BIT09}, we apply the
procedure directly to the original line integral appearing in the
left-hand side. Of course, since these two expressions are equal,
the final formula will exactly match the HBVM($k$,$s$) method,
written in a different guise.

With this premise, by considering the first expansion in
\eqref{expan},  the conservation property  reads
\begin{equation}
\label{conservation} \sum_{j=1}^{s} \gamma_j^T \int_0^1  P_j(\tau)
\nabla H(\sigma(t_0+\tau h))\mathrm{d}\tau=0,
\end{equation}
which, as is easily checked,  is certainly satisfied if we impose
the following set of orthogonality conditions
\begin{equation}
\label{orth} \gamma_j = \eta_j \int_0^1  P_j(\tau) J \nabla
H(\sigma(t_0+\tau h))\mathrm{d}\tau, \qquad j=1,\dots,s,
\end{equation}
with $\{\eta_j\}$ suitably {\em nonzero} scaling factors that will
be defined in a while.  Then, from the second relation of
\eqref{expan} we obtain, by introducing the operator
\begin{eqnarray}\label{Lf}\lefteqn{L(f;h)\sigma(t_0+ch) =}\\ \nonumber
&& \sigma(t_0)+h\sum_{j=1}^s \eta_j \int_0^c P_j(x) \mathrm{d}x \,
\int_0^1 P_j(\tau)f(\sigma(t_0+\tau h))\mathrm{d}\tau,\qquad
c\in[0,1],\end{eqnarray} that $\sigma$ is the eigenfunction of
$L(J\nabla H;h)$ relative to the eigenvalue $\lambda=1$:
\begin{equation}\label{L}\sigma = L(J\nabla H;h)\sigma.\end{equation}
\begin{defn} Equation (\ref{L}) will be called the {\em Master Functional Equation}
defining $\sigma$.\end{defn}

\begin{rem}\label{ecc} We also observe that, from (\ref{orth}) and
the first relation in (\ref{expan}), one obtains the equations
\begin{equation}\label{ecceq} \dot\sigma(t_0+c_ih) = \sum_{j=1}^s
\eta_j P_j(c_i) \int_0^1 P_j(\tau)J\nabla H(\sigma(t_0+\tau
h))\mathrm{d}\tau, \qquad i=1,\dots,s,
\end{equation}
which may be viewed as {\em extended collocation conditions}
according to \cite[Section\,2]{IT2}, where the integrals are
(exactly) replaced by discrete sums (see, e.g.,
(\ref{discr_lin})--(\ref{betai})).
\end{rem}

To practically compute $\sigma$, we set (see (\ref{Yi}) and
(\ref{expan}))
\begin{equation}
\label{y_i}
Y_i=  \sigma(t_0+c_i h) = y_0+ h\sum_{j=1}^{s} a_{ij}  \gamma_j, \qquad
i=1,\dots,s,
\end{equation}
where
\begin{equation}\label{aij}
a_{ij}=\int_{0}^{c_i} P_j(x) \mathrm{d}x, \qquad
i,j=1,\dots,s.\end{equation} Inserting \eqref{orth} into
\eqref{y_i} yields the final formulae which define the HBVMs class
based upon the basis $\{P_j\}$:
\begin{equation}
\label{hbvm_int} Y_i=y_0+h \int_0^1  \left( \sum_{j=1}^s \eta_j
a_{ij}  P_j(\tau) \right) J \nabla H(\sigma(t_0+\tau
h))\mathrm{d}\tau, \qquad i=1,\dots,s.
\end{equation}
The constants $\{\eta_j\}$ have to be chosen in order to make the
formula consistent. Problem (\ref{y_i})--(\ref{hbvm_int}) can be
actually solved, provided that all the $\{\eta_j\}$ are different
from zero,\footnote{For example, the choice $P_j(x)=x^{j-1}$,
$j=1,\dots,s$, would lead to $\eta_1=1$, and $\eta_j=0$,
$j=2,\dots,s$. This implies that, with this choice of the basis,
$\sigma$ can only be a line (i.e., $s=1$).} and the matrix
$$\pmatrix{ccc} \int_0^{c_1} P_1(x)\mathrm{d}x
&\dots &\int_0^{c_s}
P_1(x)\mathrm{d}x\\
\vdots &&\vdots\\
\int_0^{c_1} P_s(x)\mathrm{d}x &\dots &\int_0^{c_s}
P_s(x)\mathrm{d}x\endpmatrix$$ is nonsingular (which we shall
obviously assume hereafter). Indeed, such a matrix allows us, by
using (\ref{expan}), to reformulate equation (\ref{hbvm_int}) in
terms of the (unknown) fundamental stages (\ref{Yi}). Let us now
formally set
\begin{equation}\label{fy}f(y) = J\nabla H(y),\end{equation}
and report a few examples for possible choices of the basis $\{P_j(x)\}$.
\begin{enumerate}
\item In \cite{IT2} we have chosen $\{P_1(x),\dots,P_s(x)\}$ as the Newton basis.
 This has allowed us the construction of a family of methods of order 2 and 4.
\item In \cite{BIT09}, the  abscissae $\{c_0=0\}\cup\{c_i\}\cup\{\hat c_i\}$  are
disposed according to a Lobatto distribution with $k+1$ points in
$[0,1]$ and $\{P_1(x),\dots,P_s(x)\}$ is the shifted Legendre
basis in the interval $[0,1]$.\footnote{More precisely, $P_j(x)$
is here the shifted Legendre polynomial of degree $j-1$,
$j=1,\dots,s$.} Consequently, choosing in \eqref{hbvm_int} $t_0=0,
h=1$, and $f(y(\tau)) = P_j(\tau)$, the consistency condition
yields
\begin{equation}
\label{consistency}
\eta_j=\left(\int_0^1P_j^2(x)\mathrm{d}x\right)^{-1} = 2j-1,
\quad j=1,\dots,s,
\end{equation}
which is exactly the value found in \cite{BIT09}. In such a case,
it has been shown that the resulting method has order $2s$, just
the same as the generating Lobatto IIIA method (obtained for
$k=s$).

\item In a similar way, when using the Lagrange basis $\{\ell_j(x)\}$,
by setting $f(y(\tau))\equiv1$, one obtains $\eta_j=1/b_j$ with
\begin{equation}\label{bjlj}
b_j=\int_0^1 \ell_j(x)\mathrm{d} x, \qquad \ell_j(x) =
\prod_{i=1,\,i\ne j}^s \frac{x-c_i}{c_j-c_i}.
\end{equation} Consequently, formulae \eqref{hbvm_int} become
\begin{equation}
\label{hbvm_lagr} Y_i=y_0+h \int_0^1  \left( \sum_{j=1}^s
\frac{a_{ij}}{b_j}  \ell_j(\tau) \right)   J \nabla
H(\sigma(t_0+\tau h))\mathrm{d}\tau, \qquad i=1,\dots,s.
\end{equation}
Moreover, by introducing the new variables $K_i=\dot
\sigma(t_0+c_i h)$, which are  therefore related to the $Y_i$ as
$$
Y_i=y_0+h\sum_{j=1}^s a_{ij} K_j, \qquad i=1,\dots, s,
$$
system \eqref{hbvm_lagr} can be recast in the equivalent form
provided by the {\em extended collocation conditions}
(\ref{ecceq}):
\begin{equation}
\label{hbvm_hairer} K_i= \frac{1}{b_i} \int_0^1  \ell_i(\tau) J
\nabla H(\sigma(t_0+\tau h))\mathrm{d}\tau, \qquad i=1,\dots,s.
\end{equation}
\end{enumerate}
Formulae \eqref{hbvm_hairer} are the ones  E.\,Hairer proposed in
the general case, that is for any kind of Hamiltonian function.
They were called \textit{Energy Preserving variant of Collocation
Methods} (EPCMs) \cite{hairer}. The above discussion then proves
that if the integral can be substituted by a finite sum, as in the
case where $H(y)$ is a polynomial, formulae \eqref{hbvm_int}, and
consequently   \eqref{hbvm_hairer}, become a HBVM$(k,s)$, with a
suitable value of $k$.\footnote{The same clearly happens when the
integral is only {\em approximated} by a finite sum.}

For sake of completeness, we report the nonlinear system
associated with the HBVM$(k,s)$ method, in terms of the
fundamental stages $\{Y_i\}$ and the silent stages $\{\hat Y_i\}$
(see (\ref{hYi})), by using the notation (\ref{fy}). In this
context, they  represent the discrete counterpart of
\eqref{hbvm_int}, and may be directly retrieved by evaluating, for
example, the integrals in \eqref{hbvm_int} by means of the (exact)
quadrature formula introduced in \eqref{discr_lin}:{\small
\begin{eqnarray}\nonumber
Y_i  &=& y_0+  h \left[ \sum_{l=1}^s \beta_l
\left(\sum_{j=1}^s \eta_j a_{ij} P_j(c_l)\right) f(Y_l) +
\sum_{l=1}^r \hat \beta_l \left(\sum_{j=1}^s \eta_j a_{ij}
P_j(\hat c_l)\right) f(\widehat Y_l)\right]\\
\label{hbvm_sys}\\
&=& y_0+h\sum_{j=1}^s \eta_j a_{ij}\left( \sum_{l=1}^s \beta_l
P_j(c_l)f(Y_l) + \sum_{l=1}^r\hat \beta_l P_j(\hat c_l) f(\widehat Y_l)
\right),\quad i=1,\dots,s.\nonumber
\end{eqnarray}}
From the above discussion it is clear that, in the non-polynomial
case, supposing to choose the abscissae $\{\hat c_i\}$ so that the
sums in (\ref{hbvm_sys}) converge to an integral as
$r=k-s\rightarrow\infty$, the resulting formula is
\eqref{hbvm_int}.\footnote{This obvious requirement for the
abscissae will be always assumed in the sequel.} Consequently,
EPCMs may be viewed as the limit of HBVMs family, when the
Lagrange basis is considered, as the number of silent stages grows
to infinity.

The above arguments also imply that  HBVMs may be as well applied
in the non-polynomial case since, in finite precision arithmetic,
HBVMs are indistinguishable from their limit formulae
\eqref{hbvm_int}, when a sufficient number of silent stages is
introduced. The aspect of having a {\em practical} exact integral,
for $k$ large enough, was already stressed in \cite{BIT09,
IP1,IT2}.

\section{Infinity Hamiltonian Boundary Value Methods}
\label{sect:inf} As is easily argued (and emphasized in
Remark~\ref{rem:basis}), the choice of the basis along which $\dot
\sigma(t_0+\tau h)$ is expanded (see \eqref{expan}), somehow
influences the shape of the final formulae \eqref{hbvm_int}, that
is, to two different polynomial bases there may correspond two
different families of formulae.\footnote{For example, see the
method presented in subsection \ref{ex2}.} The question then
naturally arises about the best possible choice of the basis to
consider. This issue  has been  a crucial point in devising the
class of HBVMs in \cite{BIT09} and deserves a particular
attention.\footnote{The argument presented here is the analog of
the one appearing in \cite[Remark 3.1]{IT2}.}

Indeed, we recall that our final goal is to devise methods that
make the sum \eqref{conservation}, representing the line integral,
vanish. This is accomplished by the orthogonality conditions
\eqref{orth}, whose effect is to make null each term of the sum in
\eqref{conservation}. It follows that such conditions are in
general too demanding, in that they are sufficient but not
necessary to get conservativeness. In fact, the sum could in
principle vanish even in the case when  two ore more of its terms
are different from zero.  This \textit{extra constraint}  may
affect  the general properties of the conservative methods we are
interested in, and in particular their order.

This was a problem already encountered in \cite{IT2} where the
authors realized that the use of  the Newton basis didn't assure
the expected growth of the order of the resulting method  when the
degree of the polynomial $\sigma(t_0+\tau h)$ was increased. This
barrier has been definitively overcome in \cite{BIT09}, where it
was understood that the proper polynomial basis to be used by
default was that of the shifted Legendre polynomials in the
interval $[0,1]$. We emphasize that, contrary to what happens for
the Lagrange and Newton bases, the Legendre polynomials are
orthogonal and symmetric in the interval $[0,1]$ and in addition
they are \textit{abscissae-free}, that is they by no means depend
on the specific distribution of the abscissae $\{c_i\}$ adopted.
This, in turn,  implies that the {\em Master Functional Equation}
(\ref{L}) is independent of the choice of both the abscissae
$\{c_i\}$ and $\{\hat c_i\}$: the only requirement being that
(\ref{discr_lin}) holds true.\footnote{We emphasize that this is
not the case when using, for example, the Lagrange basis.}

From the above arguments, it is clear that the orthogonality
conditions (\ref{orth}), i.e., the fulfillment of the {\em Master
Functional Equation} (\ref{L}), is only a sufficient condition for
the conservation property (\ref{conservation}) to hold. Such a
condition becomes also necessary, when the basis $\{P_j\}$ is
orthogonal.

\begin{theo}\label{ortbas} Let $\{P_j\}$ be an orthogonal basis on the
interval $[0,1]$. Then, assuming $H(y)$ to be suitably
differentiable, (\ref{conservation}) implies that each term in the sum
has to vanish.
\end{theo}
\begin{proof}
Let us consider, for simplicity, the case of an orthonormal basis,
and the expansion$$g(\tau) \equiv \nabla H(\sigma(t_0+\tau h)) =
\sum_{\ell\ge 1} \rho_\ell P_\ell(\tau), \qquad \rho_\ell
=(P_\ell,g), \qquad \ell\ge1,$$ where, in general,
$$(f,g) = \int_0^1f(\tau)g(\tau) \mathrm{d}\tau.$$
Substituting into (\ref{conservation}), yields
$$\sum_{j=1}^{s} \gamma_j^T (P_j,g)=\sum_{j=1}^{s} \gamma_j^T
\left(P_j,\sum_{\ell\ge 1} \rho_\ell P_\ell\right)=\sum_{j=1}^s
\gamma_j^T\rho_j=0.$$ Since this has to hold whatever the choice
of the function $H(y)$, one concludes that
\begin{equation}\label{ortoj}\gamma_j^T\rho_j=0, \qquad j=1,\dots,s.\QED
\end{equation}
\end{proof}\medskip

\begin{rem}\label{J}
In the case where $\{P_j\}$ is the shifted Legendre basis, from
(\ref{ortoj}) one derives that $$\gamma_j = S \rho_j, \qquad
i=1,\dots,s,$$ where $S$ is any nonsingular skew-symmetric matrix.
The natural choice $S=J$ then leads to (\ref{orth}), with
$\eta_j=(P_j,P_j)^{-1}=2j-1$.
\end{rem}

The use of the Legendre basis allows the resulting methods to have
the best order and stability properties that one can expect. This
aspect is elucidated in the two theorems and the corollary below,
which represent the main result of the present work.

Although, up to now, we have maintained the treatment of HBVMs at
a general level, it is clear that, in view of the result presented
in Theorem \ref{ordine}, when the curve $\sigma(t_0+\tau h)$ is
assumed of polynomial type, we will implicitly adopt the Legendre
basis.\footnote{Actually, the term \textit{Hamiltonian Boundary
Value Method} has been coined in \cite{BIT09}, after introducing
the Legendre basis.} This important assumption will be
incorporated in the HBVM methods from now on: if needed, the use
of any other kind of basis will be explicitly stated, in order not
to create confusion.

Taking into account the consistency conditions
\eqref{consistency}, formula (\ref{hbvm_int})--(\ref{aij}) takes the form:
\begin{equation}
\label{hbvm_inf}
Y_i=y_0+h \int_0^1  \left( \sum_{j=1}^s (2j-1)
a_{ij}  P_j(\tau) \right) J \nabla H(\sigma(t_0+\tau
h))\mathrm{d}\tau, \qquad i=1,\dots,s.
\end{equation}
If the Hamiltonian $H(y)$ is a polynomial, the integral appearing
at the right-hand side is exactly computed by a quadrature
formula, thus resulting into a HBVM($k$,$s$) method with a
sufficient number of silent stages. As already stressed in the
previous section, in the non-polynomial case such formulae
represent the limit of the sequence HBVM($k$,$s$), as $k
\rightarrow \infty$.

\begin{defn}
We call the new limit formula (\ref{hbvm_inf}) an {\em Infinity
Hamiltonian Boundary Value Method} (in short, $\infty$-HBVM or
HBVM$(\infty,s)$).
\end{defn}
We emphasize that, in the non-polynomial case, \eqref{hbvm_inf}
becomes an operative method, only after that a suitable strategy
to approximate the integral is taken into account (see the next
section for additional examples). In the present case, if one
discretizes the {\em Master Functional Equation}
(\ref{Lf})--(\ref{L}), HBVM$(k,s)$ are then obtained, essentially
by extending the discrete problem (\ref{hbvm_sys}) also to the
silent stages (\ref{hYi}). In order to simplify the exposition, we
shall use (\ref{fy}) and introduce the following notation:

\begin{equation}\label{tiyi}
\{t_i\} = \{c_i\} \cup \{\hat{c}_i\}, \quad
\{\omega_i\}=\{\beta_i\}\cup\{\hat\beta_i\},\quad
y_i = \sigma(t_0+t_ih), \quad
f_i = f(\sigma(t_0+t_ih)).
\end{equation}
The discrete problem defining the HBVM$(k,s)$ then becomes, with
$\eta_j=2j-1$,
\begin{equation}\label{hbvmks}
y_i = y_0 + h\sum_{j=1}^s \eta_j \int_0^{t_i} P_j(x)\mathrm{d}x
\sum_{\ell=1}^k \omega_\ell P_j(t_\ell)f_\ell, \qquad i=1,\dots,k.
\end{equation}
We can cast the set of equations in vector form, by introducing
the vectors $\bfy = (y_1^T,\dots,y_k^T)^T$,
$e=(1,\dots,1)^T\in\RR^k$, and the matrices $\I, \P\in\RR^{k\times
s}$, with
\begin{eqnarray}\label{IDPO}
\I_{ij} &=& \int_0^{t_i} P_j(x)\mathrm{d}x,\qquad \P_{ij} = P_j(t_i),\\
\D &=& \diag(\eta_1,\dots,\eta_s), \quad
\O=\diag(\omega_1,\dots,\omega_k),\nonumber
\end{eqnarray}
as \begin{equation}\label{rk0} \bfy = e\otimes y_0 + h(\I \D
\P^T\O)\otimes I\, f(\bfy),\end{equation} with an obvious meaning
of $f(\bfy)$. Consequently, the method can be seen as a
Runge-Kutta method with the following Butcher tableau:
\begin{equation}\label{rk}
\begin{array}{c|c}\begin{array}{c} t_1\\ \vdots\\ t_k\end{array} & \I\D\P^T\O\\
 \hline                    &\omega_1\, \dots~ \omega_k
                    \end{array}\end{equation}

\begin{rem}\label{ascisse} We observe that, provided that the matrix $\D$ is
independent of the basic abscissae $\{c_i\}$ (as in the case of the Legendre
basis), the role of such abscissae and of the silent abscissae $\{\hat c_i\}$ is
interchangeable. This is not true, for example, for the Newton and Lagrange
bases.
\end{rem}

The following result then holds true.
\begin{theo}\label{ordine} Provided that the quadrature has order at least
$2s$ (i.e., it is exact for polynomials of degree at least
$2s-1$), HBVM($k$,$s$) has order $p=2 s\equiv 2\deg(\sigma)$,
whatever the choice of the abscissae $c_1,\dots,c_s$.
\end{theo}

\begin{proof} From the classical result of Butcher (see, e.g.,
\cite[Theorem\,7.4]{HNW}), the thesis follows if the simplifying
assumptions $C(s)$, $B(p)$, $p\ge 2s$, and $D(s-1)$ are satisfied.
By looking at the method (\ref{rk0})--(\ref{rk}), one has that the
first two (i.e., $C(s)$ and $B(p)$, $p\ge 2s$) are obviously
fulfilled: the former by the definition of the method, the second
by hypothesis. The proof is then completed, if we prove $D(s-1)$.
Such condition can be cast in matrix form, by introducing the
vector $\bar{e}=(1,\dots,1)^T\in\RR^{s-1}$, and the matrices
$$Q=\diag(1,\dots,s-1),\qquad D=\diag(t_1,\dots,t_k),\qquad
V=(t_i^{j-1})\in\RR^{k\times s-1},$$ (see also (\ref{IDPO})) as
$$Q V^T\O\left(\I\D\P^T\O\right) = \left(\bar{e}\,e^T -V^TD\right)\O,$$ i.e.,

\begin{equation}\label{finito}
\P\D\I^T\O V Q = \left(e\,\bar{e}^T -DV\right).
\end{equation}
Since the quadrature is exact for polynomial of degree $2s-1$. one has

\begin{eqnarray*}
\left(\I^T\O VQ\right)_{ij} &=& \left( \sum_{\ell=1}^k \omega_\ell
\int_0^{t_\ell} P_i(x)\mathrm{d}x\,(j t_\ell^{j-1}) \right) =
\left(\int_0^1 \, \int_0^t P_i(x)\mathrm{d}x
(jt^{j-1})\mathrm{d}t\right)
\\&=& \left( \delta_{i1}-\int_0^1P_i(x)x^j\mathrm{d}x\right),
\qquad i = 1,\dots,s,\quad j=1,\dots,s-1,\end{eqnarray*} where the
last equality is obtained by integrating by parts, with
$\delta_{i1}$ the Kronecker symbol. Consequently,
\begin{eqnarray*}\left(\P\D\I^T\O V Q\right)_{ij} &=& \left(1 - \sum_{\ell=1}^s
\eta_\ell P_\ell(t_i)\int_0^1 P_\ell(x) x^j\mathrm{d}x \right) =
(1-t_i^j),\\ &&\qquad\qquad\qquad i=1,\dots,k,\quad
j=1,\dots,s-1,\end{eqnarray*} that is,
(\ref{finito}), where the last equality follows from the fact that
$$\sum_{\ell=1}^s \eta_\ell P_\ell(t)\int_0^1 P_\ell(x)
x^j\mathrm{d}x = t^j, \qquad j=1,\dots,s-1.\QED$$
\end{proof}

\medskip Concerning the stability, the following result holds true.

\begin{theo}\label{stab} For all $k$ such that the quadrature
formula has order at least $2s\equiv 2\deg(\sigma)$, HBVM($k$,$s$)
is perfectly $A$-stable, whatever the choice of the abscissae
$c_1,\dots,c_s$.
\end{theo}
\begin{proof} As it has been previously observed, a HBVM$(k,s)$ is
fully characterized by the corresponding polynomial $\sigma$
which, for $k$ sufficiently large (i.e., assuming that
(\ref{discr_lin}) holds true), satisfies the {\em Master
Functional Equation} (\ref{Lf})--(\ref{L}), which is independent
of the choice of the nodes $c_1,\dots,c_s$ (since we consider the
Legendre basis). When, in place of $f(y)=J\nabla H(y)$ we put the
test equation $f(y)=\lambda y$, we have that the collocation
polynomial of the Gauss-Legendre method of order $2s$, say
$\sigma_s$, satisfies the {\em Master Functional Equation}, since
the integrands appearing in it are polynomials of degree at most
$2s-1$, so that $\sigma=\sigma_s$. The proof completes by
considering that Gauss-Legendre methods are perfectly
$A$-stable.\QED
\end{proof}
\medskip

\medskip A worthwhile consequence of Theorems~\ref{ordine} and \ref{stab}
is that one can transfer to HBVM$(\infty,s)$  all those properties
of HBVM($k$,$s$) which are satisfied starting from a given $k \ge
k_0$ on: for example, the order and stability properties.

\begin{cor}
\label{ordineinf} Whatever the choice of the abscissae
$c_1,\dots,c_s$, HBVM$(\infty,s)$ \eqref{hbvm_inf} has order $2s$
and is perfectly $A$-stable.
\end{cor}

\begin{rem}
From the result of Corollary~\ref{ordineinf}, it follows that
HBVM$(\infty,s)$ has order $2s$ and is perfectly $A$-stable for
{\em any} choice for the abscissae $c_1,\dots,c_s$. Since such
abscissae can be arbitrarily chosen, we can formally place them at
the roots of the Gauss-Legendre polynomial of degree $s$. On the
other hand, by considering that, at such abscissae, by setting
$\{\ell_i(c)\}$ and $\{b_i\}$ the corresponding Lagrange
polynomials and quadrature weights, respectively (see
(\ref{bjlj})),
$$\frac{1}{b_i} \int_0^1 P_j(x)\ell_i(x)\mathrm{d}x = \frac{1}{b_i}\sum_{r=1}^s
b_r P_j(c_r)\ell_i(c_r) = P_j(c_i), \qquad j=1,\dots,s,$$ one
obtains (with $\eta_j=2j-1$ and by using the notation (\ref{fy})):

\begin{eqnarray*}
\sigma'(t_0+c_ih) &=&\sum_{j=1}^s \eta_j
P_j(c_i)\int_0^1P_j(\tau)f(\sigma(t_0+\tau
h))\mathrm{d}\tau\\
&=&\sum_{j=1}^s \eta_j\left(\frac{1}{b_i}\int_0^1
P_j(x)\ell_i(x)\mathrm{d}x\right) \int_0^1P_j(\tau)f(\sigma(t_0+\tau
h))\mathrm{d}\tau\\
&=&\frac{1}{b_i}\int_0^1\left(\sum_{j=1}^s
\eta_jP_j(\tau)\int_0^1 P_j(x)\ell_i(x)\mathrm{d}x\right)f(\sigma(t_0+\tau
h))\mathrm{d}\tau\\
&=& \frac{1}{b_i}\int_0^1
\ell_i(\tau)f(\sigma(t_0+\tau h))\mathrm{d}\tau, \qquad i=1,\dots,s.\\
\end{eqnarray*}
Consequently,  {\em for any choice of the abscissae $\{c_i\}$},
HBVM$(\infty,s)$ provide the same polynomial $\sigma$ as the ``optimal EPCMs''
(\ref{hbvm_hairer}) of order $2s$ \cite{hairer}.\footnote{In this sense,
they are {\em equivalent}, even though they generate different discrete
problems.} Conversely, an EPCM is optimal (i.e., it has order $2s$) {\em only}
when the abscissae $c_1,\dots,c_s$ define a quadrature formula of order at least
$2s-1$, whereas different choices result in methods of lower order
\cite[Theorem\,1]{hairer}.
\end{rem}

\begin{rem}\label{alleqrem}
We also observe that, due to the choice of the shifted Legendre
polynomial basis (see (\ref{hbvm_inf}))
$$ \mbox{HBVM}(\infty,s) = \lim_{k\rightarrow\infty}
\mbox{HBVM}(k,s),$$ whatever is the choice of the fundamental
abscissae $\{c_i\}$. Consequently, for all $k$ large enough, so
that the {\em Master Functional Equation} (\ref{L}) holds true
(e.g., in the case of a polynomial Hamiltonian $H(y)$), all
HBVM$(k,s)$ provide the same polynomial $\sigma$ of degree $s$,
{\em independently of the choice of the abscissae $\{c_i\}$}.
Hence, they are {\em equivalent} to each other. This result doesn't
change in the case where $H(y)$ is not a polynomial, provided that
$H(y)$ is sufficiently differentiable. In this case, in fact, one
formally obtains, in place of the {\em Master Functional Equation
(\ref{L})}, an equation of the form
\begin{equation}\label{Lp}
\sigma_k = L(J\nabla H;h)\sigma_k + \psi_k(h),
\end{equation}
where $\psi_k(h) = O(h^{q_k-s+2})$, $q_k$ being the degree of
precision of the quadrature at the right-hand side in
(\ref{discr_lin}), so that $q_k\rightarrow\infty$ as
$k\rightarrow\infty$. From (\ref{L}) and (\ref{Lp}), one then
obtains that as $h\rightarrow0$, assuming that $f$ is Lipschitzian
with constant $\mu$, and for a suitable constant $M$ independent
of $h$:
$$\|\sigma_k - \sigma\| \le h\mu M \|\sigma_k - \sigma\| +
\|\psi_k(h)\|,$$ i.e.,
$$ \|\sigma_k - \sigma\| \le (1-h\mu
M)^{-1}\|\psi_k(h)\| = O(h^{q_k-s+2})\rightarrow 0, \qquad
k\rightarrow\infty.$$ One then concludes that, when using finite
precision arithmetic, $\sigma_k$ is indistinguishable from
$\sigma$, for all $k$ large enough.
\end{rem}

\begin{ex} As previously mentioned, for the methods studied in
\cite{BIT09}, based on a Lobatto distribution of the nodes
$\{c_0=0,c_1,\dots,c_s\}\cup\{\hat{c}_1,\dots,\hat{c}_{k-s}\}$,
one has that $\deg(\sigma)=s$, so that the order of HBVM($k$,$s$)
turns out to be $2s$, with a quadrature satisfying
$B(2k)$.\end{ex}

\begin{ex}\label{gaussex} For the same reason, when one considers a Gauss
distribution for the abscissae
$\{c_1,\dots,c_s\}\cup\{\hat{c}_1,\dots,\hat{c}_{k-s}\}$, one also
obtains a method of order $2s$ with a quadrature satisfying
$B(2k)$. This case will be further studied in
Section~\ref{sect:gauss}.
\end{ex}

\begin{rem}\label{symrem}
Finally, we also mention that, from Remark~\ref{ascisse},
HBVM($k$,$s$) are symmetric methods,\footnote{According to the
{\em time reversal symmetry condition} defined in
\cite[p.\,218]{BT}.} provided that the abscissae $\{t_i\}$ (see
(\ref{tiyi})) are symmetrically distributed (see also
\cite{BIT09}).\end{rem}

\section{Generalization of Hamiltonian BVMs}
\label{sect:gen} The approach that has allowed the construction of
methods that conserve energy functions of polynomial type is quite
general: that is, by no means it depends on the particular vector
space generating the curve $\sigma(t)$ nor on the quadrature
technique used. As was emphasized in \cite[Section\,2]{IT2}, it
solely relies on the following two ingredients: the definition of
\textit{discrete line integral} and the \textit{extended
collocation conditions} (\ref{ecceq}), which zero the line
integral (\ref{discr_lin}).

Therefore, in a more general context, this procedure can be
formalized as follows.  One first picks a curve $\sigma(t_0+\tau
h)$, $\tau \in [0,1]$, joining two points of the phase space
$y_0=\sigma(t_0)$ and $y_1=\sigma(t_0+h)$. Such a curve is assumed
to lie in a proper finite dimensional vector space
$W=\mathrm{span} \{P_1(x),\dots,P_s(x) \}$, where now $P_j(x)$,
$j=1,\dots,s$, are any linearly independent functions. Therefore
the curves $\sigma(t)$ and $\dot \sigma(t)$ will admit an
expansion in the form \eqref{expan}.

The fundamental hypothesis, for this approach to work, is that the
choice of $W$ must guarantee that the  functions $P_j(\tau)\nabla
H(\sigma(t_0+\tau h))$ appearing in \eqref{conservation} (and
hence $\sigmd(t)^T\nabla H(\sigma(t))$) be elementary integrable,
that is they are required to admit a primitive that can be
expressed in terms of elementary functions. If this is the case,
all the steps performed to obtain \eqref{hbvm_int} may be repeated
with the integral substituted by the primitive.

This represents a generalization of what done for polynomial
Hamiltonian functions  not only because the vector space $W$ may
be generated by non-poly\-nomial functions but also because the
analytic solution of the line integral may be carried out by any
available technique. Hereafter, we report a couple of examples in
the class \eqref{hbvm_inf}.\footnote{While the method in
Section~\ref{ex1} is equivalently obtainable by  applying either
\eqref{hbvm_inf} or \eqref{hbvm_int}, the same is not true for the
fourth-order method derived in Section~\ref{ex2} where, the use of
the Lagrange basis, would produce a coefficient $b_2=0$ (appearing
as a denominator in the resulting formulae \eqref{hbvm_hairer}).}

\subsection{A method of order two}
\label{ex1} We consider a separable Hamiltonian function (for
simplicity we assume $m=1$)
\begin{equation}
\label{sepH} H(q,p) = V(p) - U(q).
\end{equation}
Let  $\sigma(t)$ be the segment joining $y_0=(q_0, p_0)^T$ to
$y_1=(q_1,p_1)^T$:
$$
\sigma(t_0+\tau h) = y_0+\tau(y_1-y_0).
$$
We have $c_0=0$, $c_1=1$, and the corresponding method
\eqref{hbvm_inf} becomes:
\begin{equation}
\label{itoh_abe_sep} \pmatrix{c} \displaystyle
\frac{q_1-q_0}{h} \\[.3cm]
 \displaystyle \frac{p_1-p_0}{h}
\endpmatrix = \pmatrix{c} \displaystyle
\int_0^1 V' (p_0+\tau(p_1-p_0))\mathrm{d} \tau \\[.3cm]
 \displaystyle \int_0^1 U' (q_0+\tau(q_1-q_0))\mathrm{d} \tau
\endpmatrix =
\pmatrix{c} \displaystyle
\frac{V(p_1)-V(p_0)}{p_1-p_0} \\[.3cm]
 \displaystyle \frac{U(q_1)-U(q_0)}{q_1-q_0}
\endpmatrix.
\end{equation}
Formula \eqref{itoh_abe_sep} is one of the simplest
\textit{discrete gradient methods}  due to Itoh and Abe \cite{IA},
whose general form, for non-separable Hamiltonian functions with
one degree of freedom, reads
\begin{equation}
\label{itoh_abe} \pmatrix{c} \displaystyle
\frac{q_1-q_0}{h} \\[.3cm]
 \displaystyle \frac{p_1-p_0}{h}
\endpmatrix =  J
\pmatrix{c} \displaystyle  \frac{H(q_1,p_0)-H(q_0,p_0)}{q_1-q_0}
 \\[.3cm]
 \displaystyle \frac{H(q_1,p_1)-H(q_1,p_0)}{p_1-p_0}
\endpmatrix.
\end{equation}
The vector appearing at the right-hand side of \eqref{itoh_abe} is
obtained by replacing the partial derivatives of $H(q, p)$ with
increments along the $q$ and $p$ axes. Method \eqref{itoh_abe} is
in general first order and nonsymmetric. However, when confined to
separable Hamiltonian systems,  it turns out to be second order
and symmetric.\footnote{A generalization of \eqref{itoh_abe}
introduced in  \cite{MF} also becomes method \eqref{itoh_abe_sep}
when applied to Hamiltonian functions in the form \eqref{sepH}.}

\subsection{A method of order four}
\label{ex2} \rm To construct a method of order four in the form
\eqref{hbvm_inf} applied to \eqref{sepH},  we pick  a curve
$\sigma(t)$ of degree two, based upon the abscissae $c_0=0$,
$c_1=1/2$, and $c_2=1$. Such a method has  been already described
in \cite{IT2} for polynomial Hamiltonian functions: here we
consider its generalization to the non-polynomial case. Setting
$Y_1=(q_{1/2}, p_{1/2})^T$ and, observing that $Y_2=(q_1,p_1)^T$,
the two components of the curve $\sigma(t_0+\tau h)$ are
\begin{equation}
\label{sigma12}
\pmatrix{c} \sigma_1(t_0+\tau h) \\[.25cm]
\sigma_2(t_0+\tau h) \endpmatrix = \pmatrix{c}
2(q_0-2q_{1/2}+q_1)\tau^2 -(3q_0-4q_{1/2}+q_1) \tau +q_0
 \\[.25cm]
2(p_0-2p_{1/2}+p_1)\tau^2 -(3p_0-4p_{1/2}+p_1) \tau +p_0
 \endpmatrix.
\end{equation}
Consequently, 
\eqref{hbvm_inf} becomes
\begin{equation}
\label{order4lob}
\begin{array}{l}
Y_1 \equiv \pmatrix{c} \displaystyle
 q_{1/2} \\[.45cm]
 \displaystyle  p_{1/2}
\endpmatrix = \pmatrix{c} \displaystyle
q_0 \\[.45cm]
 \displaystyle p_0
\endpmatrix+h \pmatrix{c} \displaystyle
\int_0^1 ( -\frac{3}{2}\tau+\frac{5}{4} )V' (\sigma_2(t_0+\tau h) )\mathrm{d} \tau \\[.45cm]
 \displaystyle \int_0^1 ( -\frac{3}{2}\tau+\frac{5}{4} ) U' (\sigma_1(t_0+\tau h))\mathrm{d} \tau
\endpmatrix, \\[1cm]
Y_2 \equiv \pmatrix{c} \displaystyle
 q_{1} \\[.45cm]
 \displaystyle  p_{1}
\endpmatrix = \pmatrix{c} \displaystyle
q_0 \\[.45cm]
 \displaystyle p_0
\endpmatrix+h \pmatrix{c} \displaystyle
\int_0^1 V' (\sigma_2(t_0+\tau h) )\mathrm{d} \tau \\[.45cm]
 \displaystyle \int_0^1  U' (\sigma_1(t_0+\tau h))\mathrm{d} \tau
\endpmatrix.
\end{array}
\end{equation}
Substituting \eqref{sigma12} into \eqref{order4lob} we obtain a
system in the unknowns $q_{1/2}$, $p_{1/2}$, $q_{1}$, $p_{1}$.
Looking at \eqref{order4lob}, we realize that even in the simpler
case of a system deriving from a Hamiltonian function in the form
\eqref{sepH}, the  elementary integrability of the integrals in
\eqref{hbvm_inf} is not a priori guaranteed. This means that, in
this case, we cannot arrive at a general formula analogous to
\eqref{itoh_abe_sep}, in terms of $U(q)$ and $V(p)$.

On the other hand, in several cases of interest, such primitive
can be explicitly computed: hereafter we report a significant
example, which we shall use later in the numerical tests in
Section~\ref{sect:gauss}.

\begin{ex}\label{lotkaex} \rm The role of this example is also to
show that, when finite precision arithmetic is used, it may be
{\em not convenient} to use the {\em infinite version} of the
methods, even if the integrals can be analytically evaluated. This
will be evident from the numerical results in Section~\ref{evai}.
The system we consider is the one defined by the Hamiltonian
function
\begin{equation}
\label{lotka_ham} H(q,p) = a(\log q - q) + b(\log p - p),
\end{equation}
where $a$ and $b$ are positive constants. The associated  system
\eqref{hamilode} reads
\begin{equation}
\label{poisson} 
\displaystyle \dot q = \hphantom{-} b\,
\left(\frac{1}{p}-1\right), \qquad \displaystyle \dot p = -
a\,\left(\frac{1}{q}-1\right).
\end{equation}
This system  is strictly related to the Lotka-Volterra model
\begin{equation}
\label{lotka} \left\{ \begin{array}{l}
\dot q = \hphantom{-} b\, q \,(1-p), \\
\dot p = - a\, p\,(1-q),
\end{array}
\right.
\end{equation}
in that system \eqref{lotka} may be recast as the Poisson system
$\dot y = \frac{1}{\eta(q,p)} J \nabla H(y)$, where
$\eta(q,p)=-\frac{1}{qp}$ is called {\em integrating factor}.

Systems \eqref{lotka} and \eqref{poisson} share the same
Hamiltonian function \eqref{lotka_ham} as first integral and,
consequently, they share the same curves as trajectories in the
phase plane. Method \eqref{order4lob} applied to \eqref{lotka}
reads

\medskip
{\footnotesize
\begin{equation}
\label{Y1_lotka}
\begin{array}{l}
\pmatrix{c} \displaystyle \frac{q_{1/2}-q_0}{h/2} \\[.5cm]
\displaystyle \frac{p_{1/2}-p_0}{h/2} \endpmatrix =   \pmatrix{c}
-b+\frac{3}{4}b\frac{\log(|p_0/p_1|)}{p_0-2p_{1/2}+p_1}
+\frac{1}{2}\frac{b}{C_1}\frac{p_0-8p_{1/2}+7p_1}{p_0-2p_{1/2}+p_1}
 \\[.125cm]   \cdot \left( \mathrm{arctanh}(\frac{-3p_0+4p_{1/2}-p_1}{C_1}) -
 \mathrm{arctanh}(\frac{p_0-4p_{1/2}+3p_1}{C_1})  \right)
\\[.5cm]
a-\frac{3}{4}a\frac{\log(|q_0/q_1|)}{q_0-2q_{1/2}+q_1}
-\frac{1}{2}\frac{a}{C_2}\frac{q_0-8q_{1/2}+7q_1}{q_0-2q_{1/2}+q_1}
 \\[.125cm]   \cdot \left( \mathrm{arctanh}(\frac{-3q_0+4q_{1/2}-q_1}{C_2}) -
 \mathrm{arctanh}(\frac{q_0-4q_{1/2}+3q_1}{C_2})  \right) \endpmatrix
\end{array},
\end{equation}
} \\{\footnotesize
\begin{equation}
\label{Y2_lotka}
\begin{array}{l}
\pmatrix{c} \displaystyle \frac{q_{1}-q_0}{h} \\[.5cm]  \displaystyle
\frac{p_{1}-p_0}{h} \endpmatrix =   \pmatrix{c}
-b-\frac{2b}{C_1} \left[
\mathrm{arctanh}(\frac{p_0-4p_{1/2}+3p_1}{C_1}) +
\mathrm{arctanh}(\frac{3p_0-4p_{1/2}+p_1}{C_1}) \right]
\\[0.5cm]
a+\frac{2a}{C_2} \left[
\mathrm{arctanh}(\frac{q_0-4q_{1/2}+3q_1}{C_2}) +
\mathrm{arctanh}(\frac{3q_0-4q_{1/2}+q_1}{C_2}) \right]
\endpmatrix
\end{array},
\end{equation}
} where
$$ \left\{ \begin{array}{l}
C_1=(p_0^2+16p_{1/2}^2+p_1^2-8p_0p_{1/2}-2p_0p_1-8p_{1/2}p_1)^{1/2}, \\[.25cm]
C_2=(q_0^2+16q_{1/2}^2+q_1^2-8q_0q_{1/2}-2q_0q_1-8q_{1/2}q_1)^{1/2}.
\end{array}
\right.
$$
\end{ex}

\section{HBVMs based upon Gauss quadrature}
\label{sect:gauss}

As anticipated in Example~\ref{gaussex}, we now study the
properties of the HBVM$(k,s)$ which is defined over the set of $k$
distinct abscissae,
$$\{t_1,\dots,t_k\} \equiv \{c_1,\dots,c_s\} \cup
\{\hat{c}_1,\dots,\hat{c}_{k-s}\},$$ coinciding with the
Gauss-Legendre nodes in $[0,1]$, i.e., the roots of the shifted
Legendre polynomial of degree $k$. The corresponding polynomial
$\sigma$ has then degree $s$. By virtue of Theorems~\ref{ordine}
and \ref{stab} (see also Remark~\ref{symrem}), such methods are
symmetric, perfectly $A$-stable, and of order $2s$. They reduce to
Gauss-Legendre collocation methods, when $k=s$, and are exact for
polynomial Hamiltonian functions of degree $\nu$, provided that

\begin{equation}\label{knu}
k\ge \frac{\nu s}2.\end{equation} By recalling what stated in
Remark~\ref{alleqrem},  for all $k$ sufficiently large so that
(\ref{discr_lin}) holds, HBVM$(k,s)$ based on the $k$
Gauss-Legendre abscissae in $[0,1]$ are {\em equivalent} to
HBVM$(k,s)$ based on $k+1$ Lobatto abscissae in $[0,1]$ (see
\cite{BIT09}), since both methods define the same polynomial
$\sigma$ of degree $s$.\footnote{In the non-polynomial case, they
converge to the same HBVM$(\infty,s)$, as $k\rightarrow\infty$.}

As matter of fact, we have run HBVM$(k,s)$ based on Gauss-Legendre
nodes, and HBVM$(k,s)$ based on the Lobatto nodes, obtaining the
same results on the polynomial test problems reported in
\cite{BIT09}, which are briefly recalled in the sequel.

\subsection{Test problem 1}  Let us consider the problem
characterized by the polynomial Hamiltonian (4.1) in \cite{Faou},
\begin{equation}\label{fhp}
H(p,q) = \frac{p^3}3 -\frac{p}2 +\frac{q^6}{30} +\frac{q^4}4
-\frac{q^3}3 +\frac{1}6,
\end{equation}

\begin{figure}[hp]
\centerline{\includegraphics[width=0.7\textwidth,height=6cm]{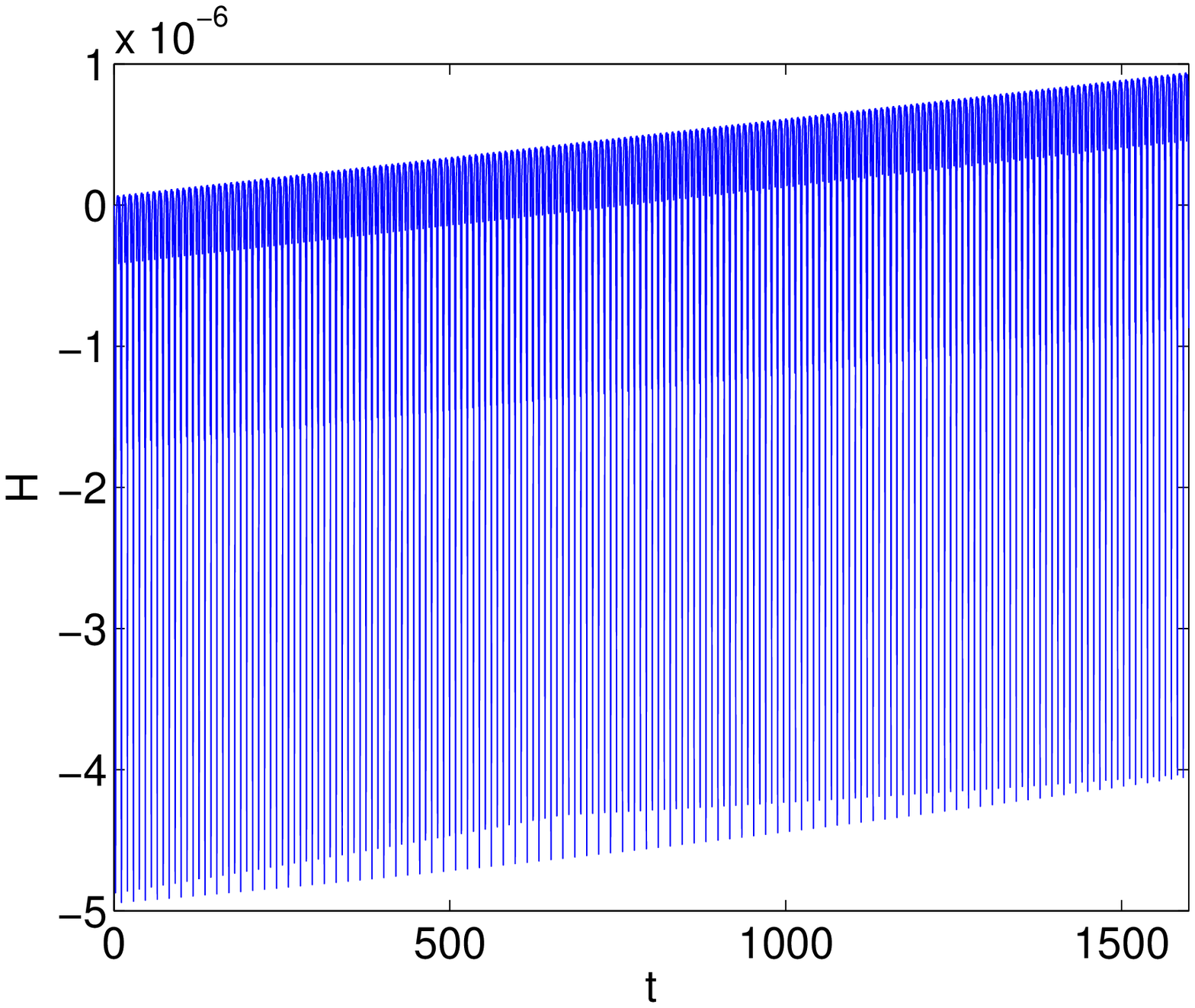}}
\caption{\protect\label{faoufig0} Fourth-order Lobatto IIIA
method, $h=0.16$, problem (\ref{fhp}): drift in the
Hamiltonian.}\medskip

\centerline{\includegraphics[width=0.7\textwidth,height=6cm]{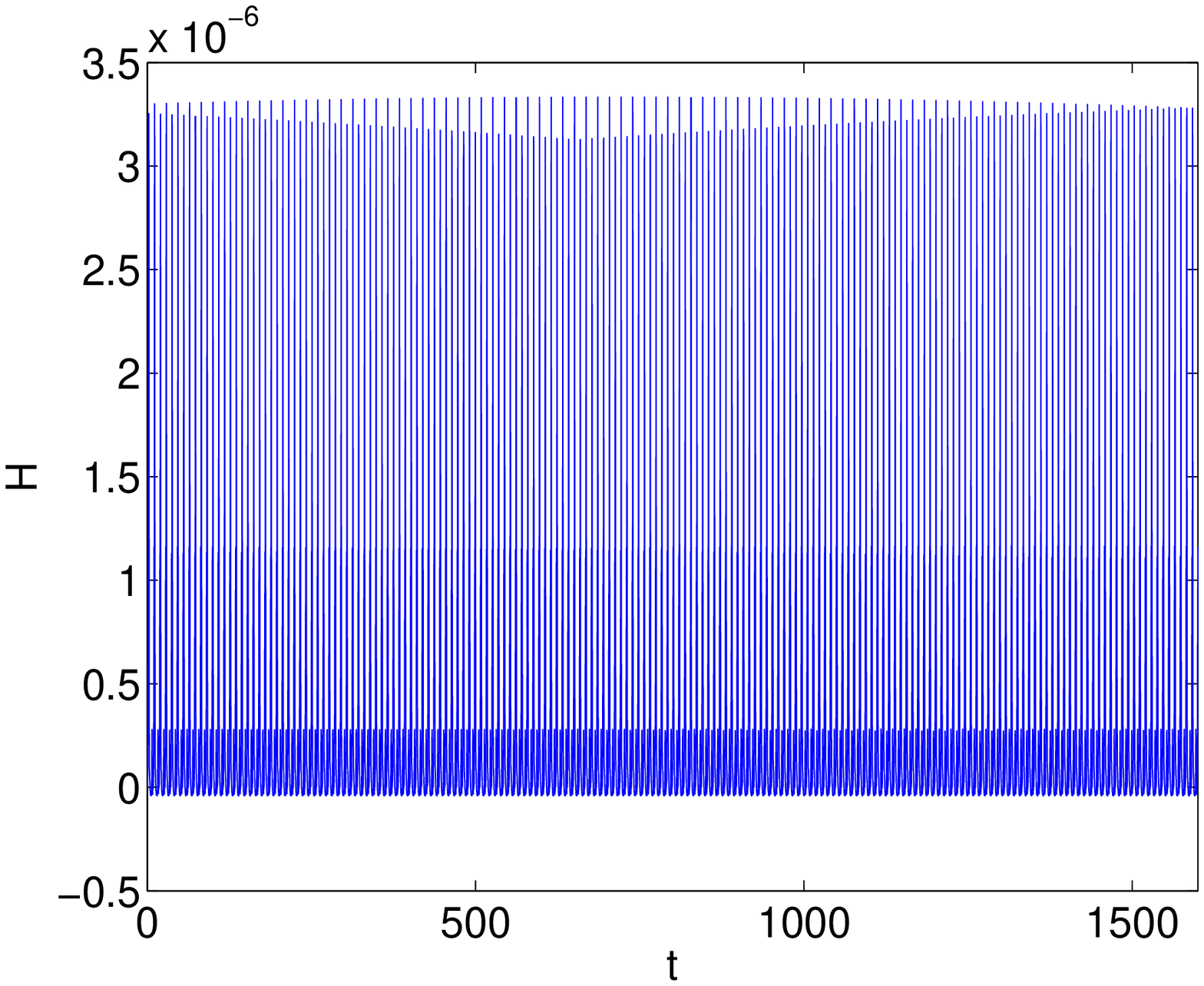}}
\caption{\protect\label{faoufig} Fourth-order Gauss-Legendre
method, $h=0.16$, problem (\ref{fhp}): $H\approx
10^{-6}$.}\medskip

\centerline{\includegraphics[width=0.7\textwidth,height=6cm]{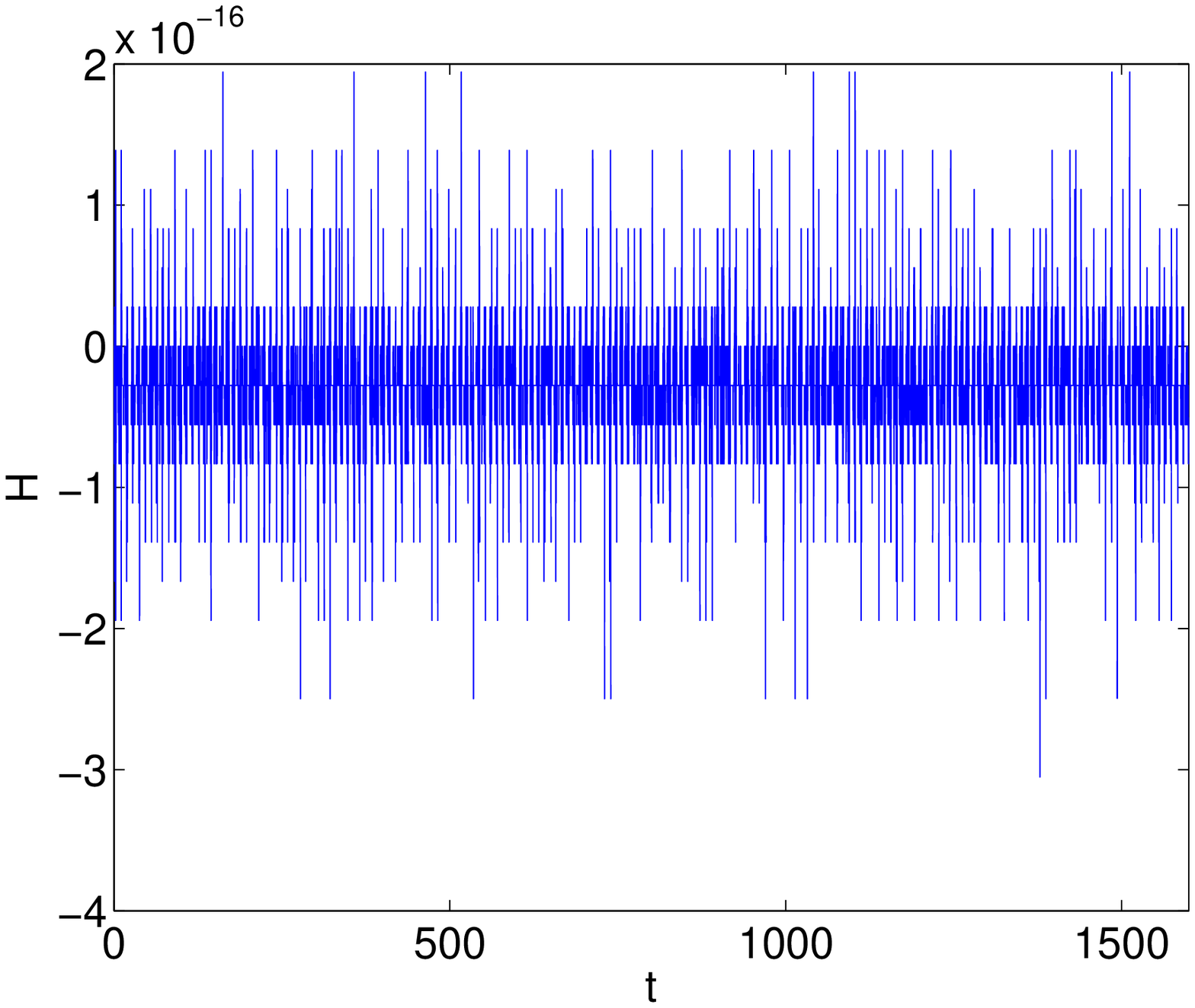}}
\caption{\protect\label{faoufig1} Fourth-order HBVM(6,2) method,
$h=0.16$, problem (\ref{fhp}): $H\approx 10^{-16}$.}
\end{figure}

\no having degree $\nu=6$, starting at the initial point
$y_0\equiv (q(0),p(0))^T=(0,1)^T$, so that $H(y_0)=0$. For such a
problem, in \cite{Faou} it has been  experienced a numerical drift
in the discrete Hamiltonian, when using the fourth-order Lobatto
IIIA method with stepsize $h=0.16$, as confirmed by the plot in
Figure~\ref{faoufig0}. When using the fourth-order Gauss-Legendre
method the drift disappears, even though the Hamiltonian is not
exactly preserved along the discrete solution, as is shown by
the plot in Figure~\ref{faoufig}. On the other hand, by using the
fourth-order HBVM(6,2) with the same stepsize, the Hamiltonian
turns out to be preserved up to machine precision, as shown in
Figure~\ref{faoufig1}, since such method exactly preserves
polynomial Hamiltonians of degree up to 6. In such a case, the
numerical solutions obtained by using the Lobatto nodes
$\{c_0=0,c_1,\dots,c_6=1\}$ or the Gauss-Legendre nodes
$\{c_1,\dots,c_6\}$ are the same.

\subsection{Test problem 2} The second test problem, having a
highly oscillating solution, is the Fermi-Pasta-Ulam problem (see
\cite[Section\,I.5.1]{hairer06gni}), defined by the Hamiltonian
\begin{equation}\label{fpu}
H(p,q) = \frac{1}2\sum_{i=1}^m\left(p_{2i-1}^2+p_{2i}^2\right)
+\frac{\omega^2}4\sum_{i=1}^m\left(q_{2i}-q_{2i-1}\right)^2
+\sum_{i=0}^m\left(q_{2i+1}-q_{2i}\right)^4,
\end{equation}

\begin{figure}[hp]
\centerline{\includegraphics[width=0.7\textwidth,height=6cm]{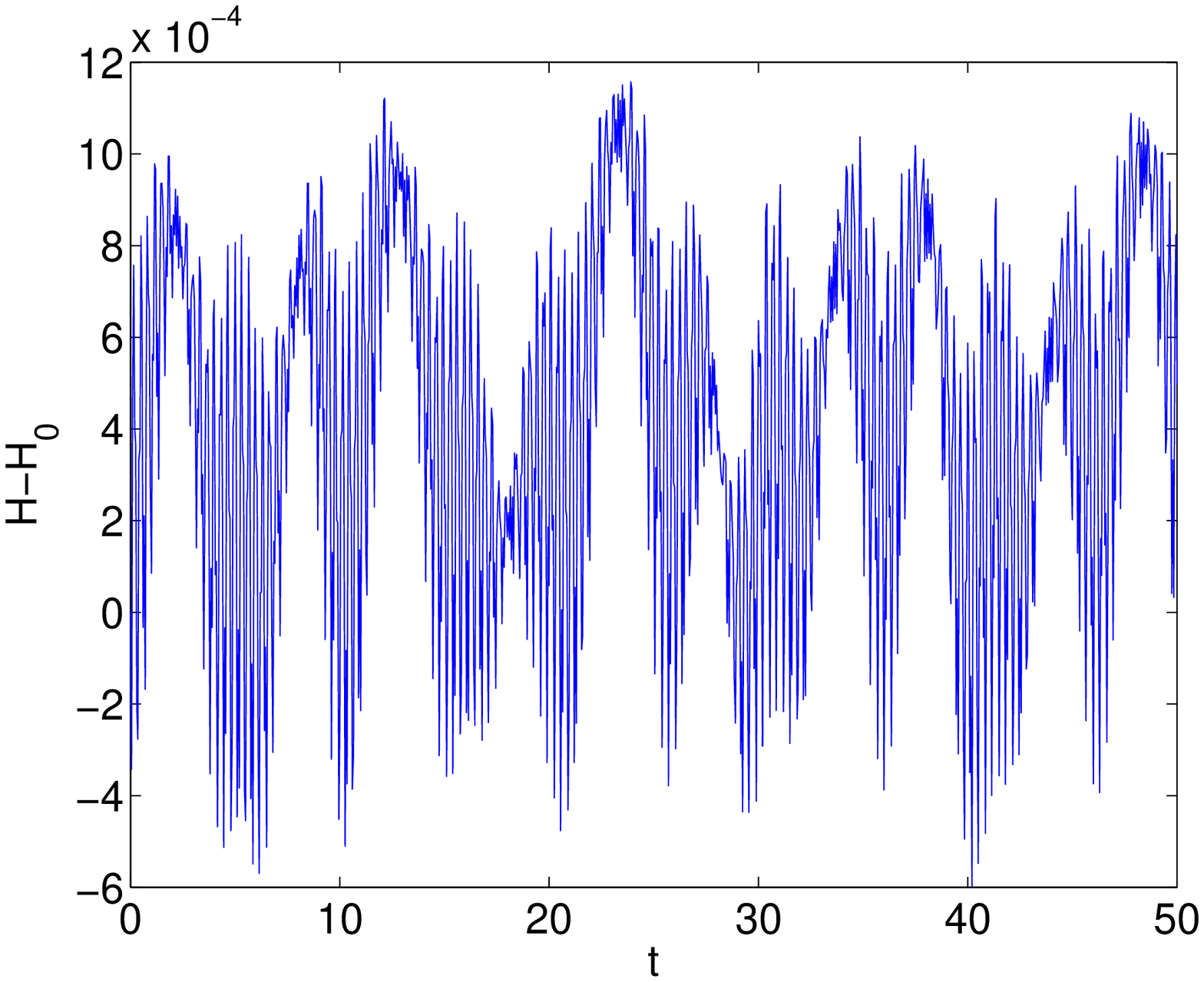}}
\caption{\protect\label{fpufig0} Fourth-order Lobatto IIIA method,
$h=0.05$, problem (\ref{fpu}): $|H-H_0|\approx 10^{-3}$.}\medskip

\centerline{\includegraphics[width=0.7\textwidth,height=6cm]{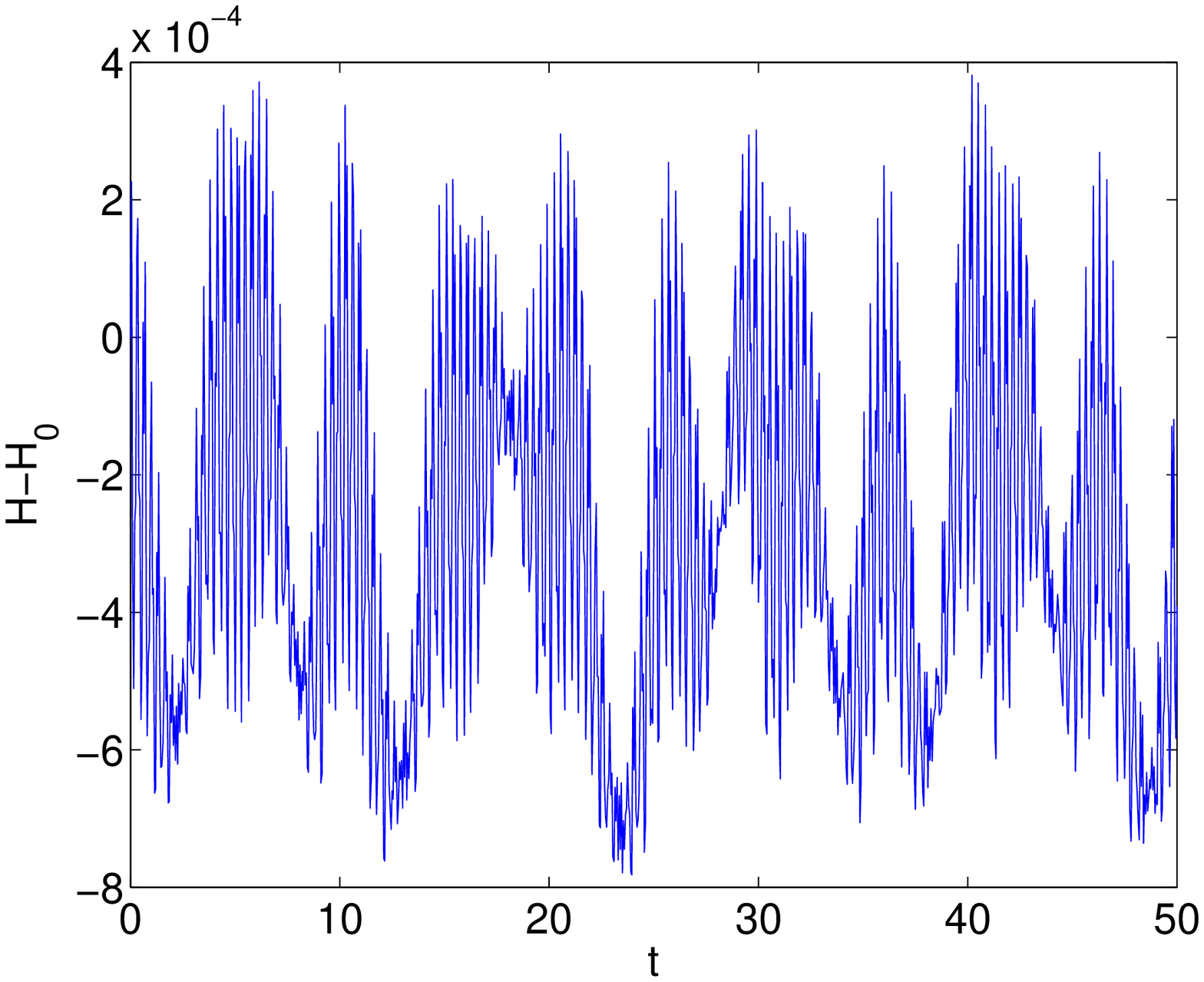}}
\caption{\protect\label{fpufig} Fourth-order Gauss-Legendre
method, $h=0.05$, problem (\ref{fpu}): $|H-H_0|\approx
10^{-3}$.}\medskip

\centerline{\includegraphics[width=0.7\textwidth,height=6cm]{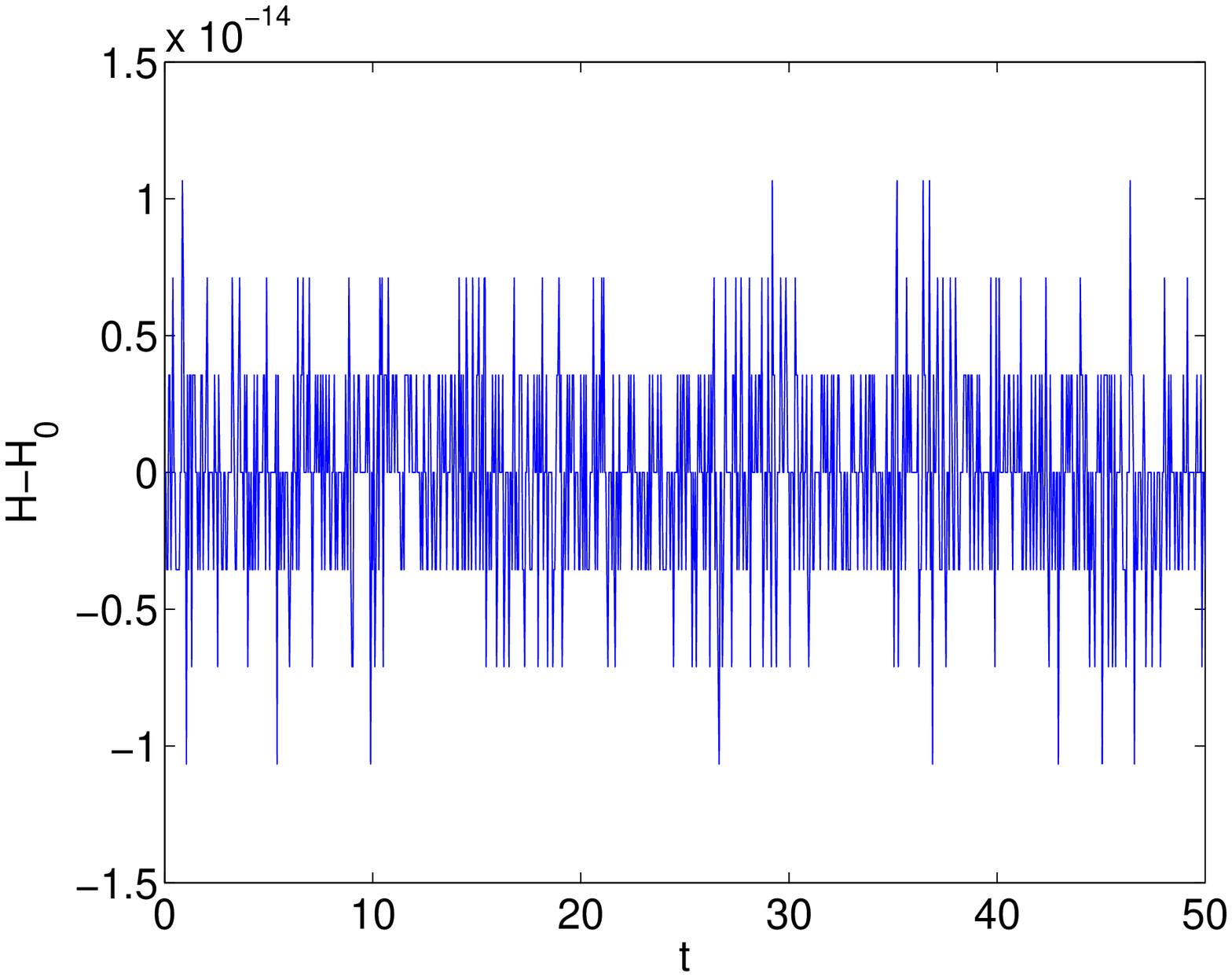}}
\caption{\protect\label{fpufig1} Fourth-order HBVM(4,2) method,
$h=0.05$, problem (\ref{fpu}): $|H-H_0|\approx 10^{-14}$.}
\end{figure}

\no with $q_0=q_{2m+1}=0$, $m=3$, $\omega=50$, and starting vector
$$p_i=0, \quad q_i = (i-1)/10, \qquad i=1,\dots,6.$$ In such a
case, the Hamiltonian function is a polynomial of degree 4, so
that the fourth-order HBVM(4,2) method, either when using the Lobatto nodes
or the Gauss-Legendre nodes, is able to exactly preserve the Hamiltonian, as
confirmed by the plot in Figure~\ref{fpufig1}, obtained with stepsize
$h=0.05$. Conversely, by using the same stepsize, both the fourth-order Lobatto
IIIA and Gauss-Legendre methods provide only an approximate conservation of the
Hamiltonian, as shown in the plots in Figures~\ref{fpufig0} and \ref{fpufig},
respectively.

\begin{table}[t]
\caption{\protect\label{tab1} Maximum difference between the
numerical solutions obtained through the fourth-order HBVM$(k,2)$
methods based on Lobatto abscissae and Gauss-Legendre abscissae
for increasing values of $k$, problem (\ref{biot}), $10^3$ steps
with stepsize $h=0.1$.}\medskip \centerline{\begin{tabular}{|r|l|}
\hline $k$  &  $h=0.1$ \\
\hline
2  & $3.97 \cdot 10^{-1}$   \\
4  & $2.29 \cdot 10^{-3}$\\
6  & $2.01 \cdot 10^{-8}$\\
8  & $1.37 \cdot 10^{-11}$\\
10 & $5.88 \cdot 10^{-13}$\\
\hline
\end{tabular}}
\end{table}

\subsection{Test problem 3 (non-polynomial Hamiltonian)} In the
previous examples, the Hamiltonian function was a polynomial.
Nevertheless, as observed above, also in this case HBVM($k$,$s$)
are expected to produce  a {\em practical} conservation of the
energy when applied to systems defined by a non-polynomial
Hamiltonian function that can be locally well approximated by a
polynomial. As an example, we consider the motion of a charged
particle in a magnetic field with Biot-Savart
potential.\footnote{\,This kind of motion causes the well known
phenomenon of {\em aurora borealis}.} It is defined by the
Hamiltonian \cite{BIT09}
\begin{eqnarray}\label{biot}
\lefteqn{H(x,y,z,\dot{x},\dot{y},\dot{z}) = }\\&&\frac{1}{2m}
\left[ \left(\dot{x}-\aa\frac{x}{\varrho^2}\right)^2 +
\left(\dot{y}-\aa\frac{y}{\varrho^2}\right)^2 +
\left(\dot{z}+\aa\log(\varrho)\right)^2\right],\nonumber
\end{eqnarray}

\no with $\varrho=\sqrt{x^2+y^2}$, $\aa= e \,B_0$,  $m$ is the
particle mass, $e$ is its charge, and $B_0$ is the magnetic field
intensity. We have used the values $$m=1, \qquad e=-1, \qquad
B_0=1,$$with starting point
$$x = 0.5, \quad y = 10, \quad z = 0, \quad
\dot{x} =  -0.1, \quad \dot{y} = -0.3, \quad \dot{z} = 0.$$

\begin{figure}[hp]
\centerline{\includegraphics[width=0.7\textwidth,height=6cm]{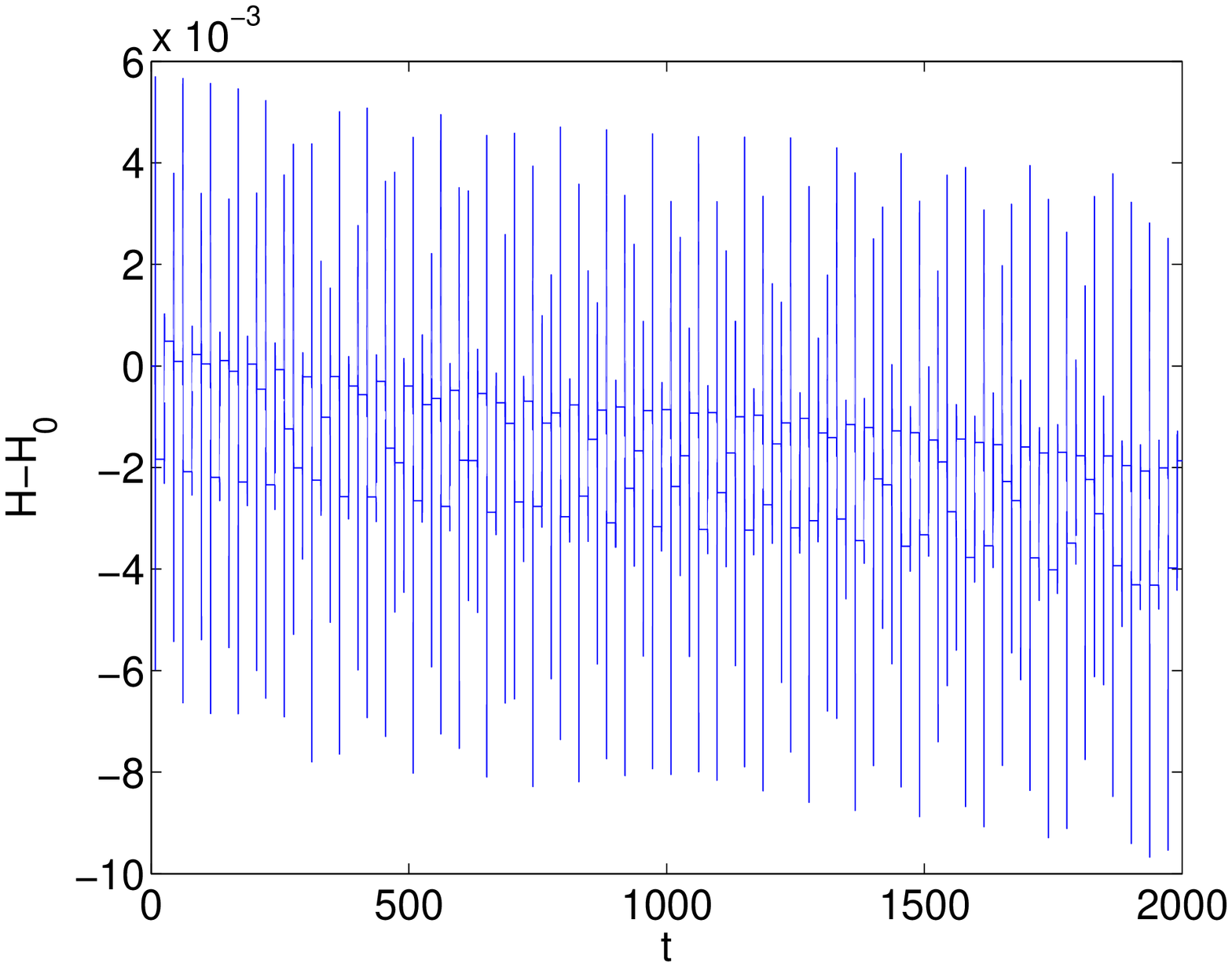}}
\caption{\protect\label{biotfig0} Fourth-order Lobatto IIIA
method, $h=0.1$, problem (\ref{biot}): drift in the
Hamiltonian.}\medskip

\centerline{\includegraphics[width=0.7\textwidth,height=6cm]{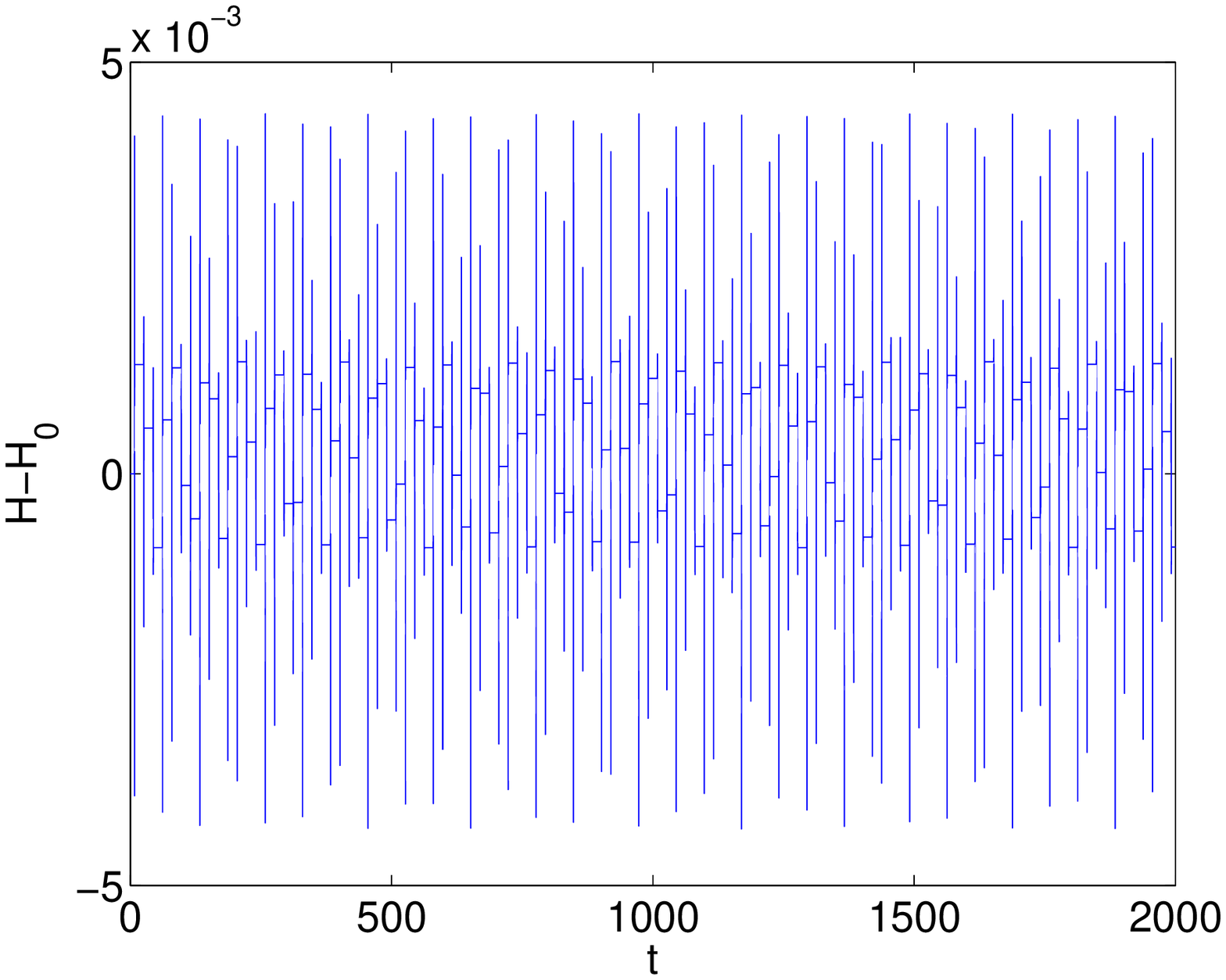}}
\caption{\protect\label{biotfig1} Fourth-order Gauss-Legendre
method, $h=0.1$, problem (\ref{biot}): $|H-H_0|\approx
10^{-3}$.}\medskip

\centerline{\includegraphics[width=0.7\textwidth,height=6cm]{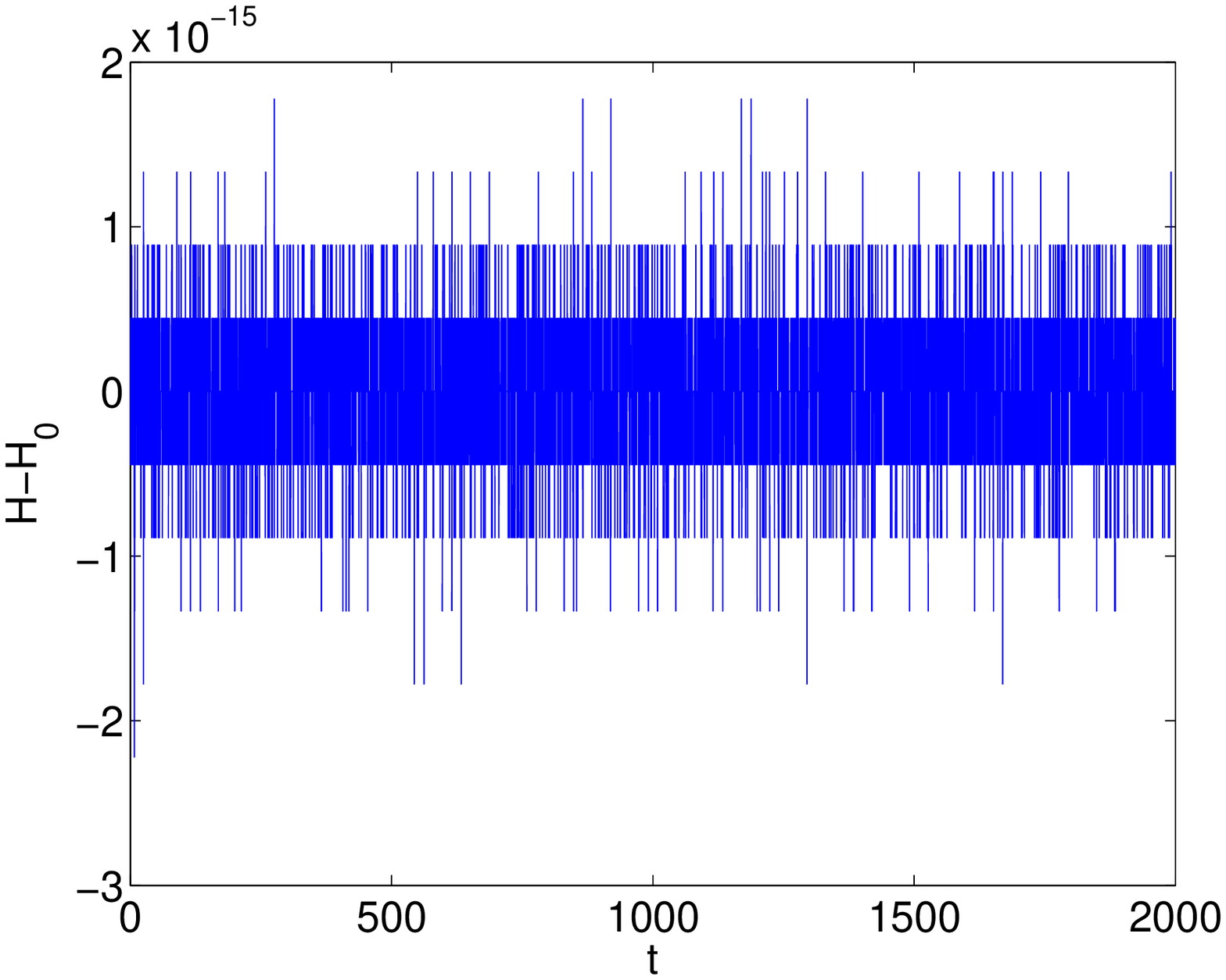}}
\caption{\protect\label{biotfig2} Fourth-order HBVM(6,2) method,
$h=0.1$, problem (\ref{biot}): $|H-H_0|\approx 10^{-15}$.}
\end{figure}

\no By using the fourth-order Lobatto IIIA method, with stepsize
$h=0.1$, a drift is again experienced in the numerical solution,
as is shown in Figure~\ref{biotfig0}. By using the fourth-order
Gauss-Legendre method with the same stepsize, the drift disappears
even though, as shown in Figure~\ref{biotfig1}, the value of the
Hamiltonian is preserved within an error of the order of
$10^{-3}$. On the other hand, when using the HBVM(6,2) method with
the same stepsize, the error in the Hamiltonian decreases to an
order of $10^{-15}$ (see Figure~\ref{biotfig2}), thus giving a
practical conservation. Finally, in Table~\ref{tab1} we list the
maximum absolute difference between the numerical solutions over
$10^3$ integration steps, computed by the HBVM$(k,2)$ methods
based on Lobatto abscissae and on Gauss-Legendre abscissae, as $k$
grows, with stepsize $h=0.1$. As expected, the difference tends to
0, as $k$ increases, since the two sequences of methods tend to
the same limit, given by the HBVM$(\infty,2)$ (see
(\ref{hbvm_inf}) with $s=2$).

\subsection{Test problem 4 (non-polynomial
Hamiltonian)}\label{evai} We finally solve the Hamiltonian system
\eqref{poisson} by using the Itho-Abe method \eqref{itoh_abe_sep},
the fourth-order formula \eqref{Y1_lotka}--\eqref{Y2_lotka}, and
the HBVM($10$,$2$), which has order four and degree of precision
$10$ (that is, according to (\ref{knu}), it precisely conserves
the energy of polynomial Hamiltonians of degree up to $10$). We
have set $a=b=1$ in formula \eqref{lotka_ham}, and integrated over
a time interval $[0, 5000]$ with  stepsize $h = 0.5$ and
$(q_0,p_0) = (0.5, 0.5)$ as initial condition.

\begin{figure}[hp]
\centerline{\includegraphics[width=0.55\textwidth,height=6cm]{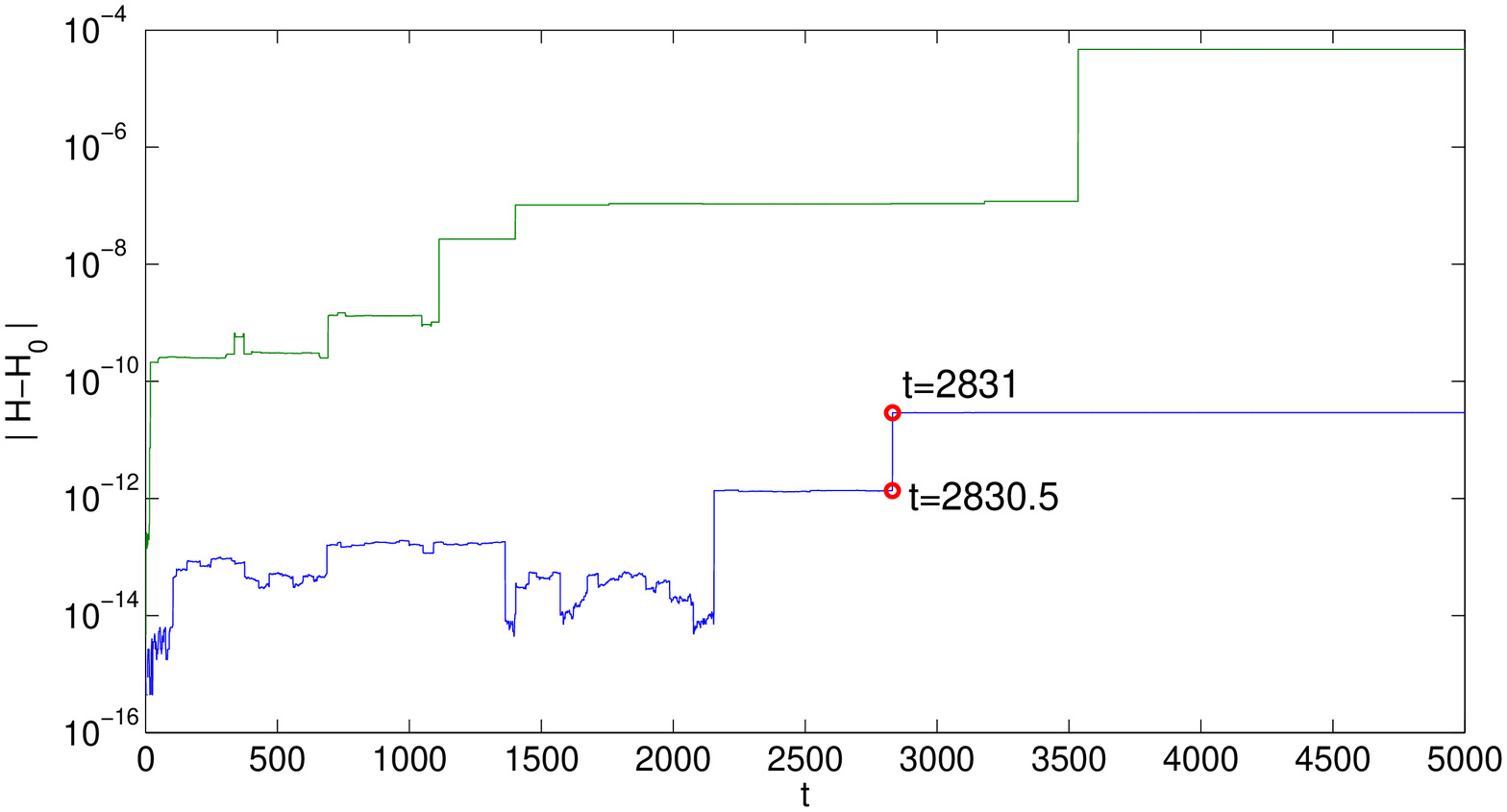}
\hspace*{-.5cm}
\includegraphics[width=0.55\textwidth,height=6cm]{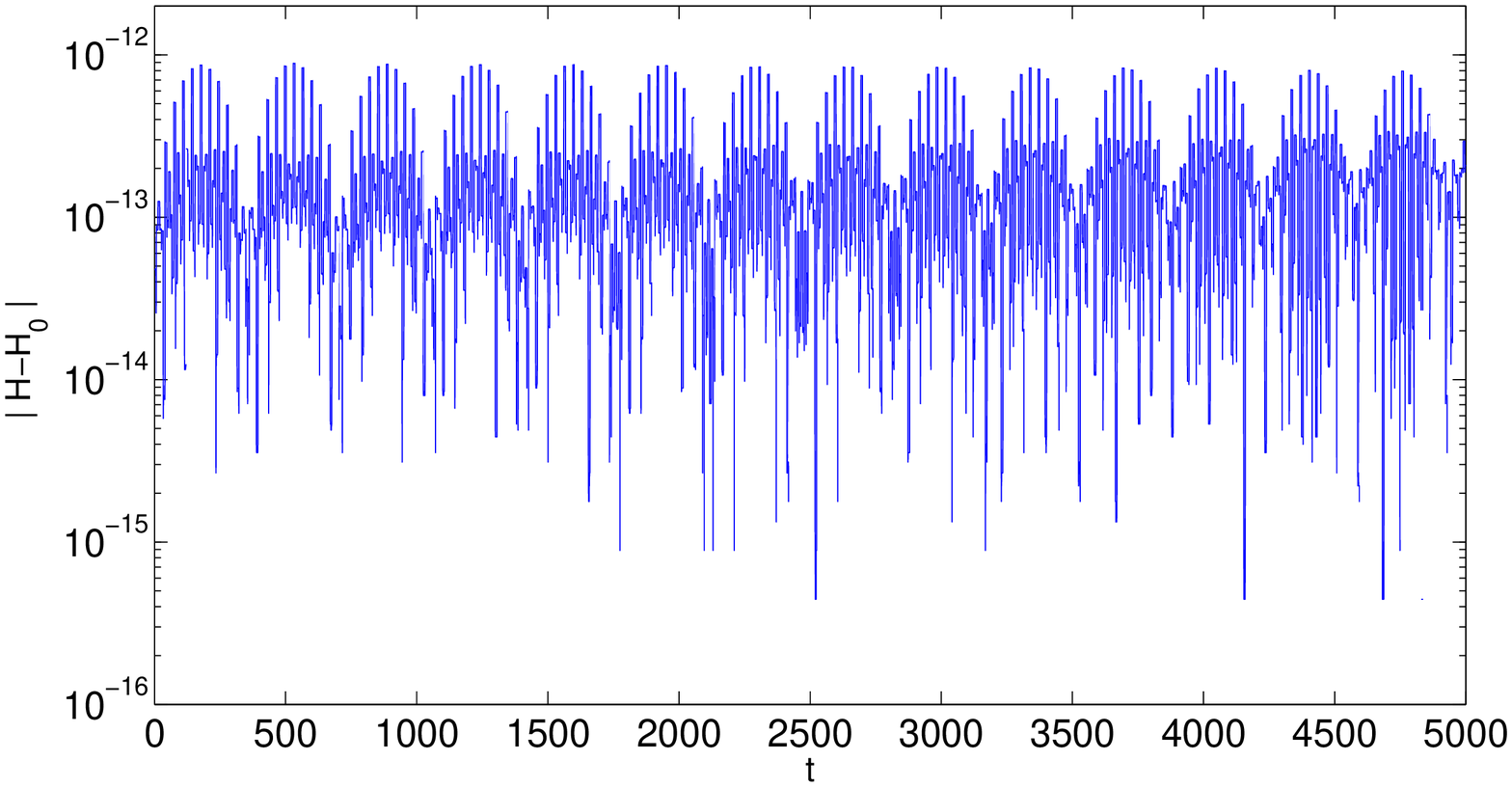}}
\caption{\protect\label{lotkafig1} Left picture: absolute error of
the Hamiltonian function \eqref{lotka_ham} evaluated along the
numerical solutions computed by the Itoh-Abe method
\eqref{itoh_abe_sep} (lower curve) and formula
\eqref{Y1_lotka}--\eqref{Y2_lotka} (upper curve). The jumps are
symptomatic of ill-conditioning of the formulae for certain values
of the solution. Right picture: the same kind of plot produced by
the HBVM formula of order $4$ and  $k=10$ Gaussian abscissae
($|H-H_0|\approx 10^{-12}$).}
\bigskip
\centerline{\includegraphics[width=1.0\textwidth]{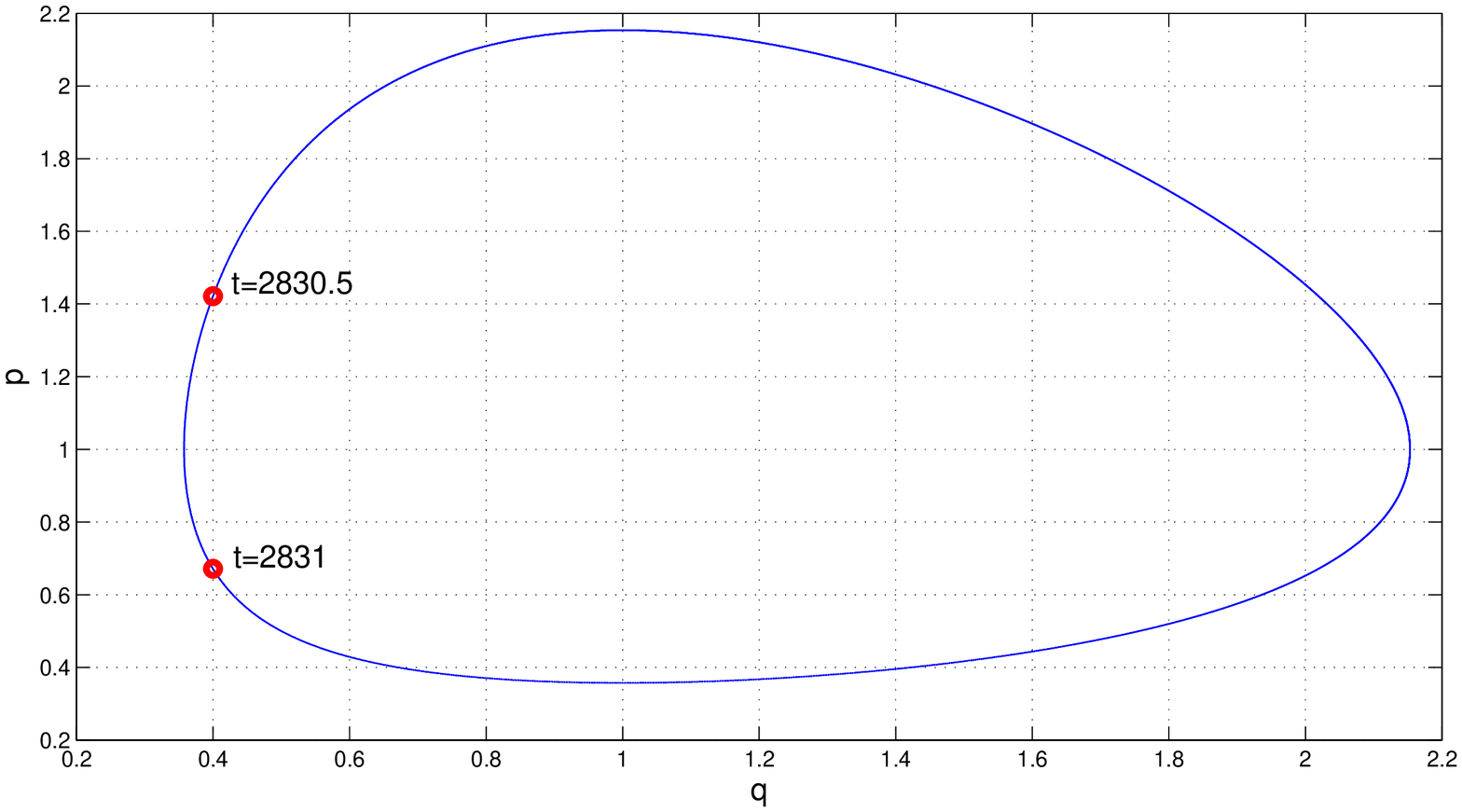}}
\caption{\protect\label{lotkafig2} Trajectory in the phase plane
computed by the Itoh-Abe method \eqref{itoh_abe_sep}. The small
circles locate the solution at the two consecutive times
$t=2830.5$ and $t=2831$. The very close values of the variable $q$
for such two points causes loss of significant digits in the
subsequent branch of the trajectory.}
\end{figure}

Figure \ref{lotkafig1} reports the numerical Hamiltonian function
associated with the three methods. The occurrence of jumps in the
first two graphs (left picture) is due to the fact that both
formulae \eqref{itoh_abe_sep} and
\eqref{Y1_lotka}--\eqref{Y2_lotka}  may become ill-conditioned for
certain values of the state vector. For example (see Figure
\ref{lotkafig2}), at the two consecutive times $t=2830.5$ and
$t=2831$, the state vectors associated with the Itoh-Abe method
\eqref{itoh_abe_sep} are, respectively,
$$
[q_1,p_1]\simeq (0.39988668,    1.4216560)^T, \qquad [q_2,p_2]\simeq (0.39988872,    0.67130503)^T,
$$
which shows that $q_1$ may be very close to $q_2$ even for large
values of $h$.  This causes some cancellation in the subtraction
at the right-hand side of \eqref{itoh_abe_sep} and, hence, a  jump
of the subsequent branch of the numerical trajectory on a
different level curve. However, since, in general, the numerical
trajectory densely fills the level curve $H(y)=H(y_0)$, it may be
argued that the occurrence of such jumps are systematic and
frequent when the dynamics is traced over a long time. The use of
finite arithmetic eventually destroys the theoretical conservation
property. A similar argument may be applied to discuss the
behavior of the fourth-order method
\eqref{Y1_lotka}--\eqref{Y2_lotka}.

Although the HBVM method does not provide a theoretical
conservation of the energy, as is the case for the above cited
methods, its behavior in finite arithmetic would suggest the
opposite (see the right picture in Figure \ref{lotkafig1}), as
already emphasized at the beginning of Example~\ref{lotkaex}.

\section{Conclusions} In this paper, the newly introduced {\em
Hamiltonian  Boundary Value Methods (HBVMs)}, a class of numerical
methods able to exactly preserve polynomial Hamiltonians of any
degree, have been re-derived in a unifying framework. Such
framework relies on the use of line integrals, which are
approximated by suitable discrete counterparts (actually, they are
exact, when the Hamiltonian is a polynomial). In this context, the
limit of the methods, as the number of the so called {\em silent
stages} tends to infinity, is easily obtained. When the underlying
polynomial basis upon which the HBVM is constructed is the
Lagrange basis, such limit formulae coincide with the recently
introduced  {\em Energy Preserving variant of Collocation
Methods}; if instead one uses the shifted Legendre polynomial
basis, the corresponding HBVMs have the highest possible order and
so do their limit formulae, the {\em Infinity Hamiltonian Boundary
Value Methods}, independently of the considered abscissae. Any
limit formula, when discretized, fall into the HBVMs class.
Possible extensions of the approach have been also sketched, along
with a number of numerical tests. Such tests confirm that, in the
limit of the silent stages tending to infinity, all HBVMs with $s$
(unknown) {\em fundamental stages} tend to the {\em same} limit
method, which is characterized by the eigenfunction (relative to
the unit eigenvalue) of a certain operator, which is independent
of the choice of the abscissae.

\section{Acknowledgements} This work emerged from fruitful
discussions with Ernst Hairer during the international conference
``ICNAAM 2009'', Rethymno, Crete, Greece, 18--22 September 2009.
It has also benefited from his many subsequent comments. We wish
also to thank an anonymous referee for the useful comments.

\end{document}